\documentclass[11pt]{article}
\usepackage{amsmath,amsthm,amssymb,color,epic,eepic, amsfonts,
cite
}
\renewcommand{\baselinestretch}{1.2}

\def\baselinestretch{1.4}
\newlength{\minitwocolumn}
\newcommand{\Z}{{\Bbb Z}} %??
\newcommand{\R}{{\Bbb R}} %??
\newcommand{\C}{{\Bbb C}} %??
\newcommand{\N}{{\Bbb N}} %???
\newcommand{\FF}{{\Bbb F}} %??
\newcommand{\PP}{{\Bbb P}} %???
\newcommand{\F}{{\mathcal F}}

\newcommand{\ta}{\tilde{a}}

\newcommand{\cD}{{\mathcal D}}
\newcommand{\cA}{{\mathcal A}}
\newcommand{\cU}{{\mathcal U}}

\newcommand{\cB}{{\mathcal B}}
\newcommand{\B}{{\mathcal B}}
\newcommand{\cI}{{\mathcal I}}

\newcommand{\cN}{{\mathcal N}}
\newcommand{\cL}{{\mathcal L}}
\newcommand{\cK}{{\mathcal K}}
\newcommand{\cR}{{\mathcal R}}
\newcommand{\cP}{{\mathcal P}}
\newcommand{\cE}{{\mathcal E}}
\newcommand{\cM}{{\mathcal M}}
\newcommand{\cO}{{\mathcal O}}
\newcommand{\cQ}{{\mathcal Q}}
\newcommand{\cS}{{\mathcal S}}

\newcommand{\cV}{{\mathcal V}}
\newcommand{\cW}{{\mathcal W}}

\newcommand{\bR}{{\overline{R}}}

\newcommand{\hL}{\widehat{L}}

\newcommand{\la}{\lambda}

\newcommand{\al}{\alpha}

\newcommand{\ep}{\epsilon}
\newcommand{\vep}{\varepsilon}

\newcommand{\tS}{\widetilde{S}}

\newcommand{\tW}{\widetilde{W}}

\newcommand{\tv}{{\widetilde{v}}}

\newcommand{\tit}{{\tilde{t}}}
\newcommand{\ti}{{\tilde{\imath}}}
\newcommand{\tj}{{\tilde{\jmath}}}
\newcommand{\tb}{{\tilde{b}}}

\newcommand{\bep}{\bar{\epsilon}}

\newcommand{\bb}{\bar{b}}
\newcommand{\bac}{\bar{c}}

\newcommand{\ha}{{\alpha}}

\newcommand{\hf}{\widehat{f}}
\newcommand{\hE}{\widehat{E}}
\newcommand{\hF}{\widehat{F}}
\newcommand{\hK}{\widehat{K}}

\newcommand{\hV}{\widehat{V}}
\newcommand{\hW}{\widehat{W}}

\newcommand{\bh}{{\bar{\h}}}

\newcommand{\hrho}{
{\rho}}

\newcommand{\Hom}{\operatorname{Hom}}

\newcommand{\nn}{{\nonumber}}
\newcommand{\bea}{\begin{eqnarray}}
\newcommand{\ena}{\end{eqnarray}}
\newcommand{\beit}{\begin{itemize}}
\newcommand{\enit}{\end{itemize}}

\newcommand{\be}{\begin{eqnarray*}}
\newcommand{\en}{\end{eqnarray*}}
\newcommand{\lb}[1]{\label{#1}}

\newcommand{\ds}[1]{{\displaystyle #1 }}

\newcommand{\End}{\mathrm{ End}}

\newcommand{\id}{\mathrm{ id}}
\newcommand{\ind}{\mathrm{ ind}}
\newcommand{\Attr}{\mathrm{ Attr}}

\newcommand{\Gr}{\mathrm{Gr}}

\newcommand{\Ell}{\mathrm{Ell}}
\newcommand{\E}{\mathrm{E}}
\newcommand{\Pic}{\mathrm{Pic}}
\newcommand{\Stab}{\mathrm{Stab}}
\newcommand{\gC}{\mathfrak{C}}
\newcommand{\pt}{\mathrm{pt}}
\newcommand{\Lie}{\mathrm{Lie}\ }

\def\infq4p#1{{(#1;q^4,p)_\infty}}

\newcommand{\tot}{\, \widetilde{\otimes}\, }

\newcommand{\mmatrix}[1]{\begin{matrix} #1 \end{matrix}}
\newcommand{\mat}[1]{\left(\mmatrix{#1}\right)}

\font\teneufm=eufm10
\font\seveneufm=eufm7
\font\fiveeufm=eufm5
\newfam\eufmfam
\textfont\eufmfam=\teneufm
\scriptfont\eufmfam=\seveneufm
\scriptscriptfont\eufmfam=\fiveeufm

\let\goth\frak
\newcommand{\slth}{\widehat{\goth{sl}}_2}
\newcommand{\slt}{\goth{sl}_2}

\newcommand{\slnh}{\widehat{\goth{sl}}_N}
\newcommand{\sln}{\goth{sl}_N}
\newcommand{\g}{\goth{g}}

\newcommand{\Bqla}{{{\mathcal B}_{q,\lambda}}}

\newcommand{\gl}{{\goth{gl}}}
\newcommand{\gln}{{\goth{gl}_N}}
\newcommand{\glnh}{\widehat{\goth{gl}}_N}

\newcommand{\h}{\goth{h}}
\newcommand{\gh}{\widehat{\goth{g}}}

\newcommand{\gS}{\goth{S}}

\makeatletter
\@addtoreset{equation}{section}
\makeatother

\newtheorem{thm}{Theorem}[section]
\newtheorem{prop}[thm]{Proposition}
\newtheorem{lem}[thm]{Lemma}
\newtheorem{cor}[thm]{Corollary}

\newtheorem{dfn}[thm]{Definition}

\begin{document}

%%%%%%%%%%%%%%%%%%%%%%%%%%%%%%%%%%%%%%%%%%%%%%%%%%%%%%%%%%%%%%%
%                                                             %
%  Elliptic Stable Envelopes and Finite-dimensional Representations of Elliptic Quantum Group           %
%                                                             %
%%%%%%%%%%%%%%%%%%%%%%%%%%%%%%%%%%%%%%%%%%%%%%%%%%%%%%%%%%%%%%%

\vspace{-1cm}
\begin{center}
{\bf\Large  Elliptic Stable Envelopes and Finite-dimensional Representations of Elliptic Quantum Group
\\[7mm] }
{\large  Hitoshi Konno}\\[6mm]
{\it  Department of Mathematics, Tokyo University of Marine Science and 
Technology, \\Etchujima, Koto-ku, Tokyo 135-8533, Japan\\
       hkonno0@kaiyodai.ac.jp}
\end{center}

\begin{abstract}
\noindent 
We construct a finite dimensional representation of the face type, i.e  dynamical, elliptic quantum group associated with  $\slnh$ 
on the Gelfand-Tsetlin basis of the tensor product of the $n$-vector representations. The result is described in a combinatorial
 way by using the partitions of $[1,n]$. We find  that the change of basis matrix from the standard  to the Gelfand-Tsetlin basis is 
 given by a specialization of the elliptic weight function obtained in the previous paper\cite{Konno17}.  Identifying the elliptic weight functions with the elliptic stable envelopes obtained by Aganagic and Okounkov, we show a correspondence of  the Gelfand-Tsetlin  bases (resp. the standard bases) to the fixed point classes  (resp.  the  stable classes)  in the equivariant elliptic cohomology $\E_T(X)$ of the cotangent bundle $X$ of the partial flag variety. As a result we obtain a geometric representation of the elliptic quantum group 
 on $\E_T(X)$.

\end{abstract}
\nopagebreak

\section{Introduction}
It has long been conjectured  that there is a parallelism between  the infinite dimensional (quantum) algebras  and (equivariant) cohomology, $\mathrm{K}$-theory,  and elliptic cohomology\cite{Gro,GKV}.  
In \cite{G91,Na94,Na98}, finite-dimensional representations of symmetrizable Kac-Moody algebras $\g$ were constructed 
in terms of homology groups of  quiver varieties. Their extension to the quantized universal enveloping algebras $U_q(\g)$ and to their affinization were constructed on equivariant $\mathrm{K}$-theory of quiver varieties\cite{GV,Gro96, Na94,Na98,Na00,Vasserot,VV99}. Note that Yangian $Y(\g)$ is obtained by replacing equivariant $\mathrm{K}$-theory by equivariant homology\cite{Vara,Na00}.  
The basic tool in these works are convolution operation and correspondences in  homology  and  equivariant $\mathrm{K}$-theory. See for example \cite{CG}.  However the elliptic case still remains conjecture. 

Stable envelopes introduced by Maulik and Okoukov\cite{MO} are  new tools to tackle this problem.  
For a quiver variety $X$, stable envelope is a  map from the equivariant cohomology of  the torus $A$-fixed point set $X^A$ to the equivariant cohomology of $X$.
It was extended to  equivariant $\mathrm{K}$-theory \cite{Okounkov} and  equivariant elliptic cohomology \cite{AO}.  
In terms of  stable envelopes  Maulik and Okoukov  constructed rational  $R$ matrices geometrically and obtained a geometric realization of the Yangian $Y_Q$ associated with a quiver $Q$\cite{MO}. Such geometric construction of $R$ matrices was  extended to the trigonometric\cite{Okounkov} and the elliptic\cite{AO} cases. 
The stable envelopes also proved to be useful in solving integrable systems\cite{MO,AO17,Smir }.   
 
This new approach was  enhanced by a discovery of a connection to the weight functions 
appearing in the hypergeometric integral solutions to the difference KZ equations. 
Gorbounov, Rim\' anyi, Tarasov and Varchenko found an identification of  rational weight functions with 
 stable envelopes for  torus-equivariant cohomology of the partial flag variety 
 $T^*\F_\la$  \cite{GRTV} and 
 extended this to the trigonometric ones for  the equivariant $\mathrm{K}$-theory \cite{RTV}. 
Furthermore  they succeeded to construct a geometric representation 
of the Yangian $Y(\gln)$ \cite{GRTV} and  the  quantum affine algebra $U_q(\glnh)$\cite{RTV} on the  equivariant 
cohomology and  the equivariant $\mathrm{K}$-theory, respectively. 
In these works,  a correspondence between  finite-dimensional representations of  quantum groups on 
the Gelfand-Tsetlin basis of the tensor product of the vector representations and  geometric representations 
is  a key to construction.  
Furthermore Felder, Rim\' anyi  and Varchenko  \cite{FRV} proposed a geometric representation of 
the dynamical elliptic quantum group $E_{\tau,y}(\gl_2)$ by using the $\slth$ type elliptic weight function 
obtained in \cite{FTV1,TV97}.

The elliptic weight functions of type $\slnh$ were derived in the previous paper\cite{Konno17} by using representation theory of 
the elliptic quantum group  $U_{q,p}(\slnh)$\cite{K98,JKOS,KK03,FKO}. The  $U_{q,p}(\slnh)$ is a  Drinfeld realization of the dynamical elliptic quantum group and is isomorphic to the central extension  of Felder's elliptic quantum group $E_{q,p}(\slnh)$\cite{Konno16}. 
Furthermore, in \cite{Konno17} their properties such as  triangularity, transition property, orthogonality, quasi-periodicity  and shuffle algebra structure  were investigated. 
Comparing these properties with those of the elliptic stable envelopes in \cite{AO}, we  conjectured that 
the elliptic weight functions can be identified with the elliptic stable envelopes. 
Some of  similar but  slightly different results  were presented in \cite{RTV17}. 

The purpose of this paper is to formulate a geometric representation of 
 the higher rank dynamical elliptic quantum group associated with $\slnh$. 
Constructing the Gelfand-Tsetlin basis  of the tensor product of the $n$-vector representations explicitly (Theorem \ref{XtW}), we obtain  finite-dimensional representations of  $E_{q,p}(\glnh)$ on it.  
In particular, we obtain an action of the half-currents of $E_{q,p}(\glnh)$ and of the associated elliptic currents of  $U_{q,p}(\slnh)$
 on the Gelfand-Tsetlin basis  (Theorem \ref{actHC} and Corollary \ref{actUqp}). The resultant representations are  described in a combinatorial way by using the partitions of $[1,n]$. 
It turns out that in the trigonometric and non-dynamical limit their combinatorial structures  coincide with those of $U_q(\slnh)$ on the equivariant $\mathrm{K}$-theory  obtained by Ginzburg and Vasserot 
\cite{GV, Vasserot} and by Nakajima\cite{Na00}.   
 
We then lift these representations to the geometric ones by  identifying the elliptic weight functions  with the 
elliptic stable envelopes. We make a direct comparison of the elliptic weight functions with the abelianization formula of the elliptic stable envelopes, which was obtained by Shenfeld \cite{Shenfeld} in the rational case  and extended to the elliptic case 
in \cite{AO}.   
We also obtain an identification of  certain specializations of the  elliptic weight functions 
  with the elliptic stable envelops restricted to the torus fixed points. 
In this restriction,  the stable envelopes play a role of the change of basis matrix elements from the stable 
classes to the fixed point classes  in $\E_T(T^*\F_\la)$.  
This allows us to define the fixed point classes  in $\E_T(T^*\F_\la)$ as transformations from the  stable classes.

 We then find  that this defining relation of the fixed point classes \eqref{IStab} is identical to 
 the change of basis relation  from the standard basis  to the Gelfand-Tsetlin basis \eqref{xiv}. 
Then a correspondence between  the Gelfand-Tsetlin bases ( resp. 
  the standard bases) and  the fixed point  classes (resp. the stable classes) in $\E_T(T^*\F_\la)$ 
yields a definition of the  actions of the half-currents of $E_{q,p}(\glnh)$ and of the  elliptic currents of  $U_{q,p}(\slnh)$ on the fixed point classes in $\E_T(T^*\F_\la)$, and provides a geometric representation of the elliptic quantum group 
 on $\E_T(T^*\F_\la)$  (Theorem \ref{geomactHC} and Corollary \ref{geomactUqp}).

In \cite{RTV17}, a similar formula for elliptic weight functions  of type $\slnh$ and their  triangularity and the orthogonality properties are presented without derivation.  There the triangular property agrees with ours but the orthogonality property seems wrong due to a lack of the dynamical shift. There are also no formulas for the shuffle algebra  in \cite{RTV17}.   
In addition, it seems that in \cite{RTV17} a different formulation of elliptic stable envelopes from the one in \cite{AO} is presented. 
The relation between them is not clear for us. However, we would like to stress that in \cite{FRV, RTV17} 
there are neither  statements on  the definition of the fixed point classes in $\E_T(T^*\F_\la)$ nor  the correspondence between the Gelfand-Tsetlin bases and the fixed point classes, which are the keys to our results. 

A part of the results  has been presented at the workshops 
``Elliptic Hypergeometric Functions in Combinatorics,
Integrable Systems and Physics'', March 20-24, 2017, ESI, Vienna,  ``Topological Field Theories, String Theory and Matrix Models",  
August 25-31, 2017, ITEP, Moscow, Infinite Analysis 17 ``Algebraic and Combinatorial Aspects in Integrable Systems",  December 4-7, 2017, Osaka City University, ``Geometric $R$-Matrices: From Geometry to Probability'', December 18-23, 2017, University of Melbourne, Creswick  and at the MSJ Autumn meeting, September 12, 2017, Yamagata University.

This paper is organized as follows. 
In Section 2 we prepare some notations including the elliptic dynamical $R$ matrices. 
We also provide defining relations  of the elliptic quantum groups $E_{q,p}(\glnh)$ and  $U_{q,p}(\glnh)$, and their basic properties.
Definition of the half-currents of $E_{q,p}(\glnh)$ 
and thier relationship to the elliptic currents of  $U_{q,p}(\glnh)$ are also exposed. 
Section 3 is devoted to a summary of the properties of the elliptic weight functions obtained in \cite{Konno17},  
such as the triangular property,  transition property, 
orthogonality, quasi-periodicity and the shuffle algebra structures.  
 In Section 4, we discuss a construction of finite-dimensional representations of $E_{q,p}(\glnh)$ and $U_{q,p}(\slnh)$ on the Gelfand-Tsetlin basis  of the tensor product of the vector representations. 
In particular we show that the change of basis matrix from the standard  to the Gelfand-Tsetlin basis is 
 given by a specialization of the elliptic weight functions. 
In Section 5, we discuss an identification between the elliptic weight functions and the elliptic stable envelopes 
and give a  geometric representation of $E_{q,p}(\glnh)$ and $U_{q,p}(\slnh)$ on $\E_T(T^*\F_\la)$. 
In Appendix A we summarize a co-algebra structure of $E_{q,p}(\glnh)$ and  $U_{q,p}(\glnh)$. 
In Appendix B we present a proof of our main Theorem \ref{actHC}. 
In Appendix C we present a direct check of Corollary \ref{actUqp} for the relation \eqref{EFFE}. 

\section{Preliminaries}
Through  this paper we follow the notations in \cite{Konno17}. 
We here list the basic ones. 
\subsection{The commutative algebra $H$}
\begin{itemize}
\item  $A=(a_{ij})\ (i,j\in \{0,1,\cdots,N-1\})$ :   the generalized Cartan matrix of  $\slnh=\widehat{\mathfrak{sl}}(N,\C)$\cite{Kac}
\item $\h=\widetilde{\h}\oplus \C d$, $\widetilde{\h}=\bar{\h}\oplus\C c$, $\bar{\h}=\oplus_{i=1}^{N-1}\C h_i$ : the Cartan subalgebra of $\slnh$  
\item $\h^*=\widetilde{\h}^*\oplus \C\delta$,  $\widetilde{\h}^*=\bh^*\oplus \C \Lambda_0$,  $\bh^*=\oplus_{i=1}^{N-1} \C \bar{\Lambda}_i$ : the dual space of $\h$ 
\item $\al_i\in \bh^*\ (1\leq i\leq N-1)$ : simple roots such that $<\al_i,h_j>=a_{ji}$  
\item  $\cQ=\oplus_{i=1}^{N-1}\Z \al_i$ :  root lattice, \quad   $\cP=\oplus_{i=1}^{N-1}\Z \bar{\Lambda}_i$ : weight lattice  
\end{itemize}
Let  $\{ \ep_j\ (1\leq j\leq N)\}$ be the orthonormal basis  in $\R^N$ 
with the inner product 
$( \ep_j, \ep_k )=\delta_{j,k}$. We set   $\bep_j=\ep_j-\sum_{k=1}^N\ep_k/N \ (1\leq j\leq N)$. Then we  have a realization   
$\al_i=\bep_i-\bep_{i+1}$ and  $\bar{\Lambda}_i=\bep_1+\cdots+\bep_i \ (1\leq i\leq N-1)$. 
For $\al\in \bar{\h}^*$ we define $h_\al\in \bar{\h}$ by $<\beta,h_\al>=(\beta,\al)$ $\forall \beta\in \bar{\h}^*$.  
We regard $\bar{\h}\oplus \bar{\h}^*$ as the Heisenberg 
algebra  by
\bea
&&~[h_{\al},\beta]
=( \al,\beta),\qquad [h_{\al},h_{\beta}]=0=
[\al,\beta]\qquad  \al, \beta \in \bar{\h}^*.\lb{HA1}
\ena
Similarly, 
\begin{itemize}
\item  $\{P_{{\alpha}}, Q_{{\beta}}\}\ 
({\alpha}, {\beta} \in \bar{\h}^*)$ : the Heisenberg algebra 
defined by 
\begin{eqnarray}
&&[P_{\al}, Q_{\beta}]=
( \al, \beta), \qquad 
[P_{\al}, P_{\beta}]=0=
[Q_{\al}, Q_{\beta}],\lb{HA2}
\end{eqnarray}
\end{itemize}
Then we define
\begin{itemize}
\item   $H=
\sum_{j=1}^{N}\C(P+h)_{\bep_j}+\sum_{j=1}^{N}\C P_{\bep_j}+\C c$, where  $(P+h)_{\bep_j}$  is an 
 abbreviation of  $P_{\bep_j}+h_{\bep_j}$.  
 \item $H^*=\widetilde{\h}^*\oplus_{j=1}^{N}\C Q_{\bep_j}
$ :  the dual space of $H$ with  paring $<\Lambda_0,c>=1
$, $<\bar{\Lambda}_i,h_j>=\delta_{i,j}$,
$<Q_\al,P_\beta>=(\al,\beta)$ and  the others vanish.

\item  $\FF=\cM_{H^*}$ : the field of meromorphic functions on $H^*$.  
\end{itemize}

\subsection{$q$-integers, infinite products and theta functions}
Let  $q$ be  generic complex numbers satisfying $|q|<1$. 
\be 
&&[n]_q=\frac{q^n-q^{-n}}{q-q^{-1}},
 \quad \\
&&(x;q)_\infty=\prod_{n=0}^\infty(1-x q^n),\quad (x;q,t)_\infty=\prod_{n,m=0}^\infty(1-x q^n t^m). 
\en
Let $r$ be a generic positive real number and  set $p=q^{2r}$. 
 We use the following Jacobi's odd theta functions.
\bea
&&[u]=q^{\frac{u^2}{r}-u}\Theta_p(z), \qquad \Theta_p(z)=(z;p)_{\infty}(p/z;p)_\infty(p;p)_\infty 
\lb{thetafull}\\
&&[u+r]=-[u], \quad [u+r\tau]=-e^{-\pi i\tau}e^{-2\pi i {u}/{r}}[u],  \quad \lb{thetaquasiperiod}
\ena
where $z=q^{2u}, p=e^{-2\pi i/\tau }$. 
For $k\in \R$, we also need the 
theta function $[u]^*$  whose elliptic nome 
is given by  $ p^*=q^{2r^*}$,  $r^*=r-k$.  We assume $r^*>0$.   
\be
&&[u]^*=q^{\frac{u^2}{r^*}-u}\Theta_{p^*}(z).  \lb{thetafulls}
\en

\subsection{The elliptic dynamical $R$-matrix of the $\slnh$ type}\lb{R-mat}

Let $\hV=\oplus_{\mu=1}^N \FF v_{\mu}$ be the $N$-dimensional vector space over $\FF$. 
The elliptic dynamical $R$-matrix  $R^\pm(z_1/z_2,\Pi)\in \End_{\C}(\hV\otimes \hV)$ of type $\slnh$ is given by
\begin{eqnarray}
R^\pm(z,\Pi)&=&\hrho^\pm(z)\bar{R}(z,\Pi),\\
\bar{R}(z,\Pi)&=&
\sum_{j=1}^{N}E_{j,j}\otimes E_{j,j}+
\sum_{1 \leq j_1< j_2 \leq N}
\left(b_{}(u,(P+h)_{j_1,j_2 })
E_{j_1,j_1}
\otimes E_{j_2,j_2}+
\bar{b}_{}(u)
E_{j_2,j_2}\otimes E_{j_1,j_1}
\right.
 \nonumber\\
&&\qquad 
\left.
+
c_{}(u,(P+h)_{j_1,j_2 })
E_{j_1,j_2}\otimes E_{j_2,j_1}+
\bar{c}_{}
(u,(P+h)_{j_1,j_2 })E_{j_2,j_1}\otimes E_{j_1,j_2}
\right),\lb{ellR}
\end{eqnarray}
where $E_{i,j}v_{\mu}=\delta_{j,\mu}v_i$, $z=q^{2u}$, $\Pi_{j,k}=q^{2(P+h)_{j,k}}$, $(P+h)_{j,k}:=(P+h)_{\bep_j}-(P+h)_{\bep_k}$, 
\bea
&&{\rho}^+(z)=q^{-\frac{N-1}{N}}z^{\frac{N-1}{rN}}
\frac{\{q^{2N}q^{-2} z\}\{q^2 z\}}{\{q^{2N} z\}\{z\}}
\frac{\{pq^{2N}/z\}\{p/z\}}{\{pq^{2N}q^{-2}/z\}\{pq^2/z\}},\lb{defrhotilde}\\[1mm]
&&{\rho}^-(z)=\rho^+(pz),\lb{defrhomtilde}\\[2mm]
&&b(u,s)=
\frac{[s+1][s-1][u]}{[s]^2[u+1]},\qquad 
\bar{b}(u)=
\frac{[u]}{[u+1]},\lb{Relements}\\
&&c(u,s)=\frac{[1][s+u]}{[s][u+1]},\qquad \bar{c}(u,s)=\frac{[1][s-u]}{[s][u+1]}\nn
\ena
and $\{z\}=(z;p,q^{2N})_\infty$.
This $R$ matrix is gauge equivalent to Jimbo-Miwa-Okado's $A_{N-1}^{(1)}$ face type Boltzmann weight\cite{JMO} 
and can be obtained \cite{Konno06} by taking the vector representation of the universal elliptic dynamical $R$ matrix\cite{JKOStg}.

The $R^\pm(z,q^{2s})$ satisfies the dynamical Yang-Baxter equation
\bea
&& R^{\pm(12)}(z_1/z_2,q^{2(s+h^{(3)})}) R^{\pm(13)}(z_1/z_3,q^{2s}) R^{\pm(23)}(z_2/z_3,q^{2(s+h^{(1)})})\nn\\
&&\qquad =R^{\pm(23)}(z_2/z_3,q^{2s}) R^{\pm(13)}(z_1/z_3,q^{2(s+h^{(2)})}) R^{\pm(12)}(z_1/z_2,q^{2s}),\lb{DYBE}
\ena
where  $q^{2h_{j,k}^{(l)}}$ acts on the $l$-th tensor space $\hV$ by  $q^{2h_{j,k}^{(l)}}  v_{\mu}=q^{2<\bep_\mu,h_{j,k}>}v_\mu$,  
and the unitarity 
\bea
&&R(z,q^{2s})R^{(21)}(z^{-1},q^{2s})=\id_{\hV\otimes \hV}. \lb{unitarity}
\ena

\subsection{The elliptic quantum groups $E_{q,p}(\glnh)$ and $U_{q,p}(\glnh)$}
We consider the dynamical elliptic quantum group realized in the two ways  $E_{q,p}(\glnh)$ and $U_{q,p}(\glnh)$. 
The elliptic algebra $E_{q,p}(\glnh)$ is a central extension of Felder's elliptic quantum group\cite{Felder1, Konno16}, whereas 
 $U_{q,p}(\glnh)$ is an elliptic and dynamical analogue \cite{K98,Konno16} of Drinfeld's new realization of 
the quantum affine algebra $U_q(\glnh)$\cite{Drinfeld}. 
For the details of the definitions we refer the reader to Sec.3 and  Appendix D.1 in \cite{Konno16}  and Appendix A in \cite{Konno17}.

\subsubsection{The elliptic algebra $E_{q,p}(\glnh)$} 
For simplicity of presentation, we treat  the elliptic algebra $E_{q,p}(\glnh)$ as a unital associative algebra 
over $\FF$ generated by ( the Laurent coefficients of) $L^+_{ij}(z)\ (1\leq i,j\leq N)$ and the central element $q^{\pm c/2}$. 
Let  $L^+(z)=\sum_{1\leq i,j\leq N} E_{ij}L^+_{ij}(z)$.   
 In  the level $k\in \R$ representation, where $c=k$,   the defining relations are given  as follows. 
\bea
&&g({P}){\hL}_{ij}(z)={\hL}_{ij}(z)\; g(P-<Q_{\bep_j},P>),\lb{rgr}\\
&&g({P+h}){\hL}_{ij}(z)={\hL}_{ij}(z)\; g(P+h-<Q_{\bep_i},P+h>),\lb{lgr}\\
&&R^{+(12)}(z_1/z_2,\Pi)\hL^{+(1)}(z_1)\hL^{+(2)}(z_2)=\hL^{+(2)}(z_2)\hL^{+(1)}(z_1)R^{+*(12)}(z_1/z_2,\Pi^*),\lb{RLL}
\ena
where $g({P+h}), g(P)\in \FF$, and 
\be
&&R^{*+}(z,\Pi^*)=\left.R^{+}(z,\Pi)\right|_{p\to p^*, r \to r^*, [u]\to [u]^*, P+h \to P}
\en
with $\Pi^*_{j,l}=q^{2P_{j,l}}$, $p^*=pq^{-2k}=q^{2r^*}$.

Setting $\hL^-(z)=\hL^+(pq^{-k}z)$, we have from Proposition D.2 in \cite{Konno16} 
\begin{eqnarray}
&&R^{-(12)}(z_1/z_2,\Pi)\widehat{L}^{-(1)}(z_1)
\widehat{L}^{-(2)}(z_2)=
\widehat{L}^{-(2)}(z_2)
\widehat{L}^{-(1)}(z_1)
R^{-*(12)}(z_1/z_2,\Pi^*),\label{RLLpmpm}\\
&&R^{\pm(12)}(q^{\pm k}z_1/z_2,\Pi)\widehat{L}^{\pm(1)}(z_1)
\widehat{L}^{\mp(2)}(z_2)=
\widehat{L}^{\mp(2)}(z_2)
\widehat{L}^{\pm(1)}(z_1)
R^{\pm*(12)}(q^{\mp k}z_1/z_2,\Pi^*),
\label{RLLpmmp}
\end{eqnarray}
where
\be
&&R^-(z,\Pi)=\rho^-(z)\bR(z,\Pi),\qquad \rho^-(z)=z^{2\frac{N-1}{N}}\rho^+(pz),\\
&&R^{*-}(z,\Pi^*)=\left.R^{-}(z,\Pi)\right|_{p\to p^*, r \to r^*, [u]\to [u]^*, P+h \to P}.
\en
 
\begin{dfn}\lb{def:Lop} 
We define the half-currents of $E_{q,p}(\glnh)$,  $F_{j,l}^\pm(z), E_{l,j}^\pm(z), (1\leq j<l \leq N)$
and $\hK_l^\pm(z)\ (1\leq l\leq N)$, as the following Gauss components of 
$\widehat{L}^\pm(z)$. 
\begin{eqnarray}
&&\widehat{L}^\pm(z)=
\left(\begin{array}{ccccc}
1&F_{1,2}^\pm(z)&F_{1,3}^\pm(z)&\cdots&F_{1,N}^\pm(z)\\
0&1&F_{2,3}^\pm(z)&\cdots&F_{2,N}^\pm(z)\\
\vdots&\ddots&\ddots&\ddots&\vdots\\
\vdots&&\ddots&1&F_{N-1,N}^\pm(z)\\
0&\cdots&\cdots&0&1
\end{array}\right)\left(
\begin{array}{cccc}
\hK^\pm_1(z)&0&\cdots&0\\
0&\hK^\pm_2(z)&&\vdots\\
\vdots&&\ddots&0\\
0&\cdots&0&\hK^\pm_{N}(z)
\end{array}
\right)\nn\\
&&\qquad\qquad\qquad\qquad\qquad\qquad\qquad\times
\left(
\begin{array}{ccccc}
1&0&\cdots&\cdots&0\\
E^\pm_{2,1}(z)&1&\ddots&&\vdots\\
E^\pm_{3,1}(z)&
E^\pm_{3,2}(z)&\ddots&\ddots&\vdots\\
\vdots&\vdots&\ddots&1&0\\
E^\pm_{N,1}(z)&E^\pm_{N,2}(z)
&\cdots&E^\pm_{N,N-1}(z)&1
\end{array}
\right).\lb{def:lhat}
\end{eqnarray}
\end{dfn}

One can express the half-currents in terms of the quantum minor determinant of the $L$ operators \cite{Konno16}.
For $1\leq a, b\leq N$, let us define $\hL^+(z)_{a,a}=(\hL^+_{i,j}(z))_{a\leq i,j\leq N}$ and
\bea
\hL^+(z)_{a,b}&=&\mat{\hL^+_{ab}(z)&\hL^+_{aa+1}(z)&\cdots&\hL^+_{aN}(z)\cr
\hL^+_{a+1b}(z)&\hL^+_{a+1a+1}(z)&\cdots&\hL^+_{a+1N}(z)\cr
\vdots&\vdots&&\vdots\cr
\hL^+_{Nb}(z)&\hL^+_{Na+1}(z)&\cdots&\hL^+_{NN}(z)\cr} \qquad \mbox{for}\ a>b
\ena
\bea
&=&\mat{\hL^+_{ab}(z)&\hL^+_{ab+1}(z)&\cdots&\hL^+_{aN}(z)\cr
\hL^+_{b+1b}(z)&\hL^+_{b+1b+1}(z)&\cdots&\hL^+_{b+1N}(z)\cr
\vdots&\vdots&&\vdots\cr
\hL^+_{Nb}(z)&\hL^+_{Nb+1}(z)&\cdots&\hL^+_{NN}(z)\cr} \qquad \mbox{for}\ a<b.
\ena
Then we obtain
\begin{thm}\lb{HCqminor}\cite{Konno16}
\be
&&\hK^\pm_j(z)={\cN'}^{-1}_{N-j+1}{q\mbox{-}\det \hL^\pm(z)_{j,j}}\left({q\mbox{-}\det \hL^\pm(zq^{-2})_{j+1,j+1}}\right)^{-1},\\
&&E^\pm_{k,j}(z)=\left({q\mbox{-}\det \hL^\pm(z)_{k,k}}\right)^{-1}{q\mbox{-}\det \hL^\pm(z)_{k,j}},\\
&&F^\pm_{j,k}(z)={q\mbox{-}\det \hL^\pm(z)_{j,k}}\left({q\mbox{-}\det \hL^\pm(z)_{k,k}}\right)^{-1} \qquad (1\leq j<k\leq N),
\en
where 
\be
&&\cN_k\rq{}=\frac{\cN_k}{\cN_{k-1}},\quad \cN_k=\prod_{1\leq a<b\leq k}\sqrt{\frac{\rho^*_0[a]^*[1]}{\rho_0[a][1]^*}},\\
&& \rho_0=-\lim_{z\to q^{-2}}\rho^+(z)\frac{[1]}{[u+1]}=q^{\frac{r-N+2}{rN}}\frac{(q^{2N};q^{2N})_\infty}{(p;p)_\infty}\frac{\{q^{2N}q^{-4}\}\{pq^2\}^2}{\{pq^4\}\{q^{2N}q^{-2}\}^2},\\
&&\rho_0^*=\rho_0|_{p\mapsto p^*, r\mapsto r^*}
\en
\end{thm}

\begin{cor}
Let us set 
\be
&&\hK(z)=\hK^+_1(z)\hK^+_2(zq^{-2})\cdots \hK^+_N(zq^{-2(N-1)}).
\en
Then the  $q$-determinant of $\hL^+(z)$ is given by 
\be
&&q\mbox{-}\det \hL^+(z)=\cN_N \hK(z)
\en
and  belongs to the center of $E_{q,p}(\glnh)$. 
\end{cor}
Moreover from Proposition 6.4 in \cite{Konno16}, we have
\begin{cor}
For $l=1,\cdots,N$, the $q$-principal minor determinant $q\mbox{-}\det\hL^+(z)_{ll}$ is  given by 
\be
&&q\mbox{-}\det\hL^+(z)_{ll}=\cN_{N-l+1}\hK^+_l(z)\hK^+_{l+1}(zq^{-2})\cdots \hK^+_N(zq^{-2(N-l)})
\en
and belongs to the center of  the subalgebra $E_{q,p}(\widehat{\gl}_{N-l+1})$ of $E_{q,p}(\glnh)$. 
\end{cor}
We define the elliptic algebra $E_{q,p}(\slnh)$ as  the quotient algebra  $E_{q,p}(\glnh)/<\hK(z)-1>$.

\subsubsection{The elliptic algebra $U_{q,p}(\slnh)$}\lb{def:Uqp}
For simplicity of presentation, we treat  the elliptic algebra $U_{q,p}(\glnh)$ as a unital associative algebra 
over $\FF$ generated by  (the  Laurent coefficients of) the elliptic currents $E_j(z), F_j(z), K^+_l(z)\ (1\leq j\leq N-1,1\leq l \leq N)$.   
In the level-$k$ $(k\in \R)$ representation, the defining relations are given in the sense of analytic continuation as follows. 
For $g(P), g(P+h)\in \FF$, 
\bea
&&g({P+h})E_j(z)=E_j(z)g({P+h}),\quad g({P})E_j(z)=E_j(z)g(P-<Q_{\al_j},P>),\lb{gegl}\\
&&g({P+h})F_j(z)=F_j(z)g(P+h-<{\al_j},P+h>),\quad g({P})F_j(z)=F_j(z)g(P),\lb{gfgl}\\
&&g({P})K^+_l(z)=K^+_l(z)g(P-<Q_{\bep_l},P>),\quad
g({P+h})K^+_l(z)=K^+_l(z)g(P+h-<Q_{\bep_l},P>),\nn\\
&&\lb{gkgl}\\
&&K^+_{ l}(z_1)K^+_{ l}(z_2)=
\frac{{\rho}^{+*}(z_1/z_2)}{{\rho}^+(z_1/z_2)}K^+_{ l}(z_2)K^+_{ l}(z_1),\lb{KlKl}\\
&&K^+_{j}(z_1)K^+_{l}(z_2)=
\frac{{\rho}^{+*}(z_1/z_2)}{{\rho}^+(z_1/z_2)}\frac{[u_1-u_2-1]^*[u_1-u_2]}{[u_1-u_2]^*[u_1-u_2-1]}K^+_{l}(z_2)K^+_{j}(z_1)
\quad (1\leq j < l\leq N),\nn\\
&&\lb{KjKl}
\ena
\bea
&&K^+_{j}(z_1)E_j(z_2)=
\frac{\left[u_1-u_2+\frac{j-N-k+1}{2}\right]^*}{\left[u_1-u_2+\frac{j-N-k+1}{2}-1\right]^*}E_j(z_2)K^+_{j}(z_1)
,
\\
&&K^+_{j+1}(z_1)E_{j}(z_2)=
\frac{\left[u_1-u_2+\frac{j-N-k+1}{2}\right]^*}{\left[u_1-u_2+\frac{j-N-k+1}{2}+1\right]^*}E_{j}(z_2)K^+_{j+1}(z_1)
,\\
&&K^+_l(z_1)E_j(z_2)=E_j(z_2)K^+_l(z_1)\qquad (l\not=j,j+1),
\\
&&K^+_{j}(z_1)F_j(z_2)=
\frac{\left[u_1-u_2+\frac{j-N+1}{2}-1\right]}{\left[u_1-u_2+\frac{j-N+1}{2}\right]}F_j(z_2)K^+_{j}(z_1)
,\\
&&K^+_{j+1}(z_1)F_{j}(z_2)=\frac{\left[u_1-u_2+\frac{j-N+1}{2}+1\right]}{\left[u_1-u_2+\frac{j-N+1}{2}\right]}F_{j}(z_2)K^+_{j+1}(z_1)
,\\
&&K^+_l(z_1)F_j(z_2)=F_j(z_2)K^+_l(z_1)\qquad (l\not=j,j+1),
\ena
\bea
&& \left[u-v-\frac{a_{ij}}{2}\right]^* E_i(z)E_j(w)
=\left[u-v+\frac{a_{ij}}{2}\right]^* 
E_j(w)E_i(z),
\lb{u7}\\
&& \left[u-v+\frac{a_{ij}}{2} \right] F_i(z)F_j(v)
=\left[u-v-\frac{a_{ij}}{2}\right]
F_j(w)F_i(z),
\lb{u8}\\
&&[E_i(z),F_j(w)]=\frac{\delta_{i,j}}{q-q^{-1}}
\left(\delta\bigl(q^{-k}\frac{z}{w}\bigr)H^-_{i}(q^{k/2}w)
-\delta\bigl(q^{k}\frac{z}{w}\bigr)H^+_{i}(q^{-k/2}w)
\right),\lb{EFFE}
\\
&&K^-_l(z)=K^+_l(pq^{-k}z),\lb{def:Km}\\
&&H^\pm_j(z)=\varrho K^{\pm}_j(q^{N-j-1}q^{k/2}z)K^{\pm}_{j+1}(q^{N-j-1}q^{k/2}z)^{-1}
\ena
\bea
&&z_1^{-\frac{1}{r^*}}
\frac{(p^*q^2z_2/z_1;p^*)_\infty}
{(p^*q^{-2}z_2/z_1;p^*)_\infty
}\left\{
({z_2}/{z})^{\frac{1}{r^*}}
\frac{(p^*q^{-1}z/z_1;p^*)_\infty 
(p^*q^{-1}z/z_2;p^*)_\infty}
{(p^*qz/z_1;p^*)_\infty 
(p^*qz/z_2;p^*)_\infty}E_i(z_1)E_i(z_2)E_j(z)\right.\nonumber\\
&&
-\left.[2]_q\frac{(p^*q^{-1}z/z_1;p^*)_\infty 
(p^*q^{-1}z_2/z;p^*)_\infty}
{(p^*qz/z_1;p^*)_\infty 
(p^*qz_2/z;p^*)_\infty}
E_i(z_{1})E_j(z)E_i(z_{2})
\right.\nn
\\
&&+\left.
(z/z_1)^{\frac{1}{r^*}}\frac{(p^*q^{-1}z_1/z;p^*)_\infty 
(p^*q^{-1}z_2/z;p^*)_\infty}
{(p^*qz_1/z;p^*)_\infty 
(p^*qz_2/z;p^*)_\infty}
E_j(z)E_i(z_{1})E_i(z_{2})
\right\}+(z_1 \leftrightarrow z_2)=0,\nn\\
&&\label{e20}
\\
&&z_1^{\frac{1}{r}}
\frac{(pq^{-2}z_2/z_1;p)_\infty}
{(pq^{2}z_2/z_1;p)_\infty
}\left\{
(z/z_2)^{\frac{1}{r}}
\frac{(pq z/z_1;p)_\infty 
(pq z/z_2;p)_\infty}
{(pq^{-1} z/z_1;p)_\infty 
(pq^{-1} z/z_2;p)_\infty}
F_i(z_1)F_i(z_2)F_j(z)\right.\nonumber\\
&&
-\left.[2]_q\frac{(pq z/z_1;p)_\infty 
(pq z_2/z;p)_\infty}
{(pq^{-1} z/z_1;p)_\infty 
(pq^{-1} z_2/z;p)_\infty}
F_i(z_{1})F_j(z)F_i(z_{2})
\right.\nn
\\
&&+\left.
(z_1/z)^{\frac{1}{r}}
\frac{(pq z_1/z;p)_\infty 
(pq z_2/z;p)_\infty}
{(pq^{-1} z_1/z;p)_\infty 
(pq^{-1} z_2/z;p)_\infty}
F_j(z)F_i(z_{1})
F_i(z_{2})
\right\} + (z_1 \leftrightarrow z_2)=0
\quad (|i-j|=1).\nonumber\\
&&~~~~~~~~~~~~~~~~~\label{e21}
\ena
where
\be
&&\varrho=\frac{(p;p)_\infty(p^*q^2;p^*)_\infty}{(p^*;p^*)_\infty(pq^2;p)_\infty}.\lb{def:kappa}
\en

\begin{prop}\cite{Konno16}
The following product  belongs to the center of $U_{q,p}(\glnh)$.
\be
&&K(z)=K^+_1(z)K^+_2(zq^{-2})\cdots K^+_N(zq^{-2(N-1)}).  
\en
\end{prop}
We define the elliptic algebra $U_{q,p}(\slnh)$ as  the quotient algebra  $U_{q,p}(\glnh)/<K(z)-1>$.

\subsubsection{Isomorphism between $E_{q,p}(\glnh)$ and $U_{q,p}(\glnh)$
}
\begin{prop}\lb{EcHc}\cite{Konno16}
{Set }
\be
&& E_j(zq^{j-N+1- c/2}):=\mu^* \left( E^+_{j+1,j}(zq^{c/2})-E^-_{j+1,j}(zq^{- c/2})\right), \\
&&F_j(zq^{j-N+1-c/2}):=\mu\left( F^+_{j,j+1}(zq^{-c/2})-F^-_{j,j+1}(zq^{c/2})\right), 
\en
where $\mu$ and $\mu^*$ satisfy 
\bea
&&\mu\mu^*=-\frac{\varrho}{q-q^{-1}}\frac{[0]'}{[1]}.\lb{mumus}
\ena 
Then identifying  $\hK^+_l(z)$ with $K^+_l(z)$, $\hK^+_l(z), E_j(z), F_j(z)$ satisfy the defining relations of $U_{q,p}(\glnh)$  in Sec.\ref{def:Uqp}. 
\end{prop}

Furthermore let us set 
\be
&&H^\pm_j(z):=\varrho \hK^\pm_j(q^{N-j-1}q^{c/2}z)\hK^\pm_{j+1}(q^{N-j-1}q^{c/2}z)^{-1}. 
\en
\begin{cor}
Under the constraint $\hK(z)=1$, the generating functions  $H^\pm_j(z), E_j(z), F_j(z)\ (1\leq j\leq N-1)$ satisfy the defining relations of $U_{q,p}(\slnh)$.  
\end{cor}

\begin{thm}\cite{Konno16}
\be
&&U_{q,p}(\glnh)\cong E_{q,p}(\glnh).
\en
\end{thm}

\begin{dfn}\lb{def:GTsubalg}
The Gelfand-Tsetlin subalgebra $\mathfrak{G}$ of $E_{q,p}(\glnh)$ as well as of $U_{q,p}(\glnh)$ is defined to be a unital subalgebra generated by 
( the Laurent  coefficients of)  $\hK^+_l(z)\ (1\leq l\leq N)$.  
\end{dfn}
 From \eqref{KlKl}-\eqref{KjKl} we obtain
\begin{prop}
The Gelfand-Tsetlin subalgebra $\mathfrak{G}$ becomes a commutative subalgebra at the level 0 $(c=0)$. 
\end{prop}

\subsubsection{Dynamical $L$ operators and half-currents}\lb{dynL}
For later convenience ( see Sec.\ref{finRep})  we introduce  the dynamical $L$ operators\cite{JKOS,KK03} by 
\be
&&L^\pm(z,P)=\hL^\pm(z)e^{\sum_{j=1}^N\pi(h_{\epsilon_j})Q_{\bep_j}},
\en
where $\pi(h_{\epsilon_j})=E_{jj}$. 
Then $L^\pm(z,P)$ commutes with the elements in $\FF$ and satisfy the full dynamical $RLL$ relations\cite{Felder1}
\bea
&&\hspace{-1cm}R^{\pm(12)}(z_1/z_2,\Pi){L}^{\pm(1)}(z_1,P)
{L}^{\pm(2)}(z_2,P+h^{(1)})=
{L}^{\pm(2)}(z_2,P)
{L}^{\pm(1)}(z_1,P+h^{(2)})
R^{\pm*(12)}(z_1/z_2,\Pi^*),\nn\\
&&\label{DRLLpmpm}\\
&&\hspace{-1.5cm}R^{\pm(12)}(q^{\pm c}z_1/z_2,\Pi)
{L}^{\pm(1)}(z_1,P)
{L}^{\mp(2)}(z_2,P+h^{(1)})=
{L}^{\mp(2)}(z_2,P)
{L}^{\pm(1)}(z_1,P+h^{(2)})
R^{\pm*(12)}(q^{\mp c}z_1/z_2,\Pi^*),\nn\\
&&\label{DRLLpmmp}
\ena
Accordingly we define the dynamical half-currents $\cK^+_l(z), E^+_{j+1,j}(z,P), F^+_{j,j+1}(z,P)\ (1\leq l \leq N, 1\leq j\leq N-1)$ 
as the corresponding Gauss coordinates of $L^\pm(z,P)$. 
Then the relation between the half-currents from $\hL^\pm(z)$ and those from $L^\pm(z,P)$ is given as follows.
\begin{prop}{\cite{Konno16}}\lb{dynamicalHC}
\be
&&\hK^\pm_l(z)=\cK^\pm_l(z)e^{-Q_{\bep_j}},\\
&&\hE^\pm_{j+1,j}(z)=e^{Q_{\bep_{j+1}}}E^\pm_{j+1,j}(z,P)e^{-Q_{\bep_j}},\\
&&\hF^\pm_{j,j+1}(z)=F^\pm_{j,j+1}(z,P).
\en
\end{prop}

\section{Elliptic Weight Functions}\lb{EWF}
In this section  we summarize  some basic properties of  the elliptic weight functions obtained in \cite{Konno17}. 

\subsection{Combinatorial notations}\lb{combNot}
Let $\ds{\hV=\bigoplus_{\mu=1}^N \FF v_\mu}$ be the same as in Sec. \ref{R-mat} and  consider its tensor product
$\hV^{\tot n}$. 
The standard basis of $\hV^{\tot n}$ is given by 
$\{v_{\mu_1}\tot \cdots \tot v_{\mu_n}\ |\ \mu_1,\cdots,\mu_n\in \{1,\cdots,N\} \}$, where $\tot$ is defined in Appendix \ref{HHopfalgebroid}.  

Let $[1,n]=\{1,\cdots,n\}$. 
For a vector $v_{\mu_1}\tot \cdots \tot v_{\mu_n}$, we define the index set $I_l:=\{\ i\in [1,n]\ |\ \mu_i=l\}$ $(l=1,\cdots,N)$ and set  $\la_l:=|I_l|$, 
$\la:=(\la_1,\cdots,\la_N)$.  
 Then $I=(I_1,\cdots,I_N)$ is a partition of $[1,n]$, i.e.  
 \be
I_1\cup \cdots \cup I_N=[1,n],\quad I_k\cap I_l=\emptyset\quad  \mbox{$(k\not=l)$}.
\en
We often denote thus obtained partition $I$ by $I_{\mu_1,\cdots\mu_n}$. 
We also write $v_I=v_{\mu_1}\tot \cdots \tot v_{\mu_n}$. 
Let $\N=\{m\in \Z|\ m\geq 0 \}$. 
For $\la=(\la_1,\cdots,\la_N)\in \N^N$ satisfying  $|\la|=\la_1+\cdots+\la_N=n$, let $\cI_\la$ be the set of all partitions  $I=(I_1,\cdots,I_N)$ of $[1,n]$  satisfying $|I_l|=\la_l\ (l=1,\cdots,N)$.
 We also set $\la^{(l)}:=\la_1+\cdots+\la_l$,  $I^{(l)}:=I_1\cup\cdots \cup I_l$ and let $I^{(l)}=:\{i^{(l)}_1< \cdots<i^{(l)}_{\la^{(l)}}\}$. 
For $I\in \cI_\la$, all vectors $v_I$ have the same weight  $\sum_{j=1}^n\bep_{\mu_j}$, which we call
 the weight associated with $\la$.  
For each $i^{(l)}_a$ $(l=1,\cdots,N,\ a=1,\cdots,\la^{(l)})$, we consider the variables  $t^{(l)}_a\equiv t(i^{(l)}_a)$ with $t^{(N)}_a=z_a$ $(a=1,\cdots,n)$, and set $t=(t^{(l)}_a)\ (l=1,\cdots,N,\ a=1,\cdots,\la^{(l)})$.

For $\la=(\la_1,\cdots,\la_N)\in \N^N, |\la|=n$, we consider in Sec.\ref{GeomRep}  
 the partial flag variety  $\F_\la=\F(\la^{(1)},\cdots,\la^{(N-1)},n)$  consisting of 
$0=\cV_0\subset \cV_1\subset \cdots\subset \cV_N=\C^n$
with $\dim_\C \cV_l=\la^{(l)}$.  A representation theoretical meaning to this parametrization is given in \cite{Konno17}.

\subsection{The elliptic weight functions  of type $\sln$}
We consider the following elliptic weight functions\cite{Konno17}.  
\bea
&&
\widetilde{W}_I(t,z,\Pi)= {\rm Sym}_{t^{(1)}}\cdots {\rm Sym}_{t^{(N-1)}}\widetilde{U}_I(t,z,\Pi),\lb{def:Wtilde}\\
&&\widetilde{U}_I(t,z,\Pi)
=\prod_{l=1}^{N-1}\prod_{a=1}^{\la^{(l)}}\left(\frac{[v^{(l+1)}_b-v^{(l)}_a+(P+h)_{\mu_s,l+1}-C_{\mu_s,l+1}(s)][1]}{[v^{(l+1)}_b-v^{(l)}_a+{1}][(P+h)_{\mu_s,l+1}-C_{\mu_s,l+1}(s)]}\right|_{i^{(l+1)}_b=i^{(l)}_a=s}\nn\\[2mm]
&&\left.\qquad\qquad\qquad\qquad\times \prod_{b=1\atop i^{(l+1)}_b>i^{(l)}_a}^{\la^{(l+1)}}\frac{[v^{(l+1)}_b-v^{(l)}_a]}{[v^{(l+1)}_b-v^{(l)}_a+{1}]}\prod_{b=a+1}^{\la^{(l)}}\frac{[v^{(l)}_a-v^{(l)}_b-1]}{[v^{(l)}_a-v^{(l)}_b]}\right),
\ena
where  we set $t^{(l)}_a=q^{2v^{(l)}_a}\ (l=1,\cdots, N-1, a=1,\cdots,\la^{(l)})$, $z_k=q^{2u_k}\ (k=1,\cdots,n)$, $v_s^{(N)}=u_s$ $(s=1,\cdots,n)$ and $C_{\mu_s,l+1}(s):=\sum_{j=s+1}^n<\bep_{\mu_j},h_{\mu_s,l+1}>$\ $(\mu_s\leq l)$. The symbol ${\rm Sym}_{t^{(l)}}$ denotes the symmetrization over the variables $t^{(l)}_1,\cdots,t^{(l)}_{\la^{(l)}}$. 

For  $I=(I_1,\cdots,I_N)\in \cI_\la$, 
let  $I_k=\{i_{k,1}<\cdots<i_{k,\la_k}\}$ $(k=1,\cdots,N)$.
Then  $C_{\mu_s,l+1}$  has the following combinatorial expression. 
\begin{prop}\lb{combC}
\be
&&
C_{\mu_s,l+1}(s)
=\left\{\mmatrix{\la_{\mu_s}-\la_{l+1}-\tilde{s}+m_{\mu_s,l+1}(s)-1&\quad \mbox{if }\ s\leq i_{l+1,\la_{l+1}}\cr
\la_{\mu_s}-\tilde{s}\qquad\qquad &\quad \mbox{if }\ s>i_{l+1,\la_{l+1}}\cr}\right.
\en
where for $s\in [1,n]$ we define $\tilde{s}$ by $i_{\mu_s,\tilde{s}}=s$  and $m_{\mu_s,l+1}(s)$  by  
\be
&&m_{\mu_s,l+1}(s)=\mbox{min}\{1\leq j\leq \la_{l+1}\ | \ s<i_{l+1,j}\ \} \quad \mbox{for }\ s \leq i_{l+1,\la_{l+1}}.
\en
\end{prop}

\subsection{Entire function version}
Let us set 
\bea
&&H_\la(t,z):=\prod_{l=1}^{N-1}\prod_{a=1}^{\la^{(l)}}\prod_{b=1}^{\la^{(l+1)}}\left[v^{(l+1)}_b-v^{(l)}_a+{1}\right].\lb{Efunc}
\ena
The following gives an entire function version of the elliptic weight function. 
\bea
&&W_I(t,z,\Pi)=H_\la(t,z)\widetilde{W}_I(t,z,\Pi)
= {\rm Sym}_{t^{(1)}}\cdots {\rm Sym}_{t^{(N-1)}}
{U}_I(t,z,\Pi),\lb{entireW}
\ena
where
\bea
&&{U}_I(t,z,\Pi)= 
\prod_{l=1}^{N-1}\prod_{a=1}^{\la^{(l)}}\left(\left.
\frac{\left[v^{(l+1)}_b-v^{(l)}_a+(P+h)_{\mu_s,l+1}-C_{\mu_s,l+1}(s)\right][1]}{[(P+h)_{\mu_s,l+1}-C_{\mu_s,l+1}(s)]}
\right|_{i^{(l+1)}_b=i^{(l)}_a=s}\right.\nn\\
&&\left.\qquad\qquad\times \prod_{b=1\atop i^{(l+1)}_b>i^{(l)}_a}^{\la^{(l+1)}}{\left[v^{(l+1)}_b-v^{(l)}_a\right]}
\prod_{b=1\atop i^{(l+1)}_b<i^{(l)}_a}^{\la^{(l+1)}}{\left[v^{(l+1)}_b-v^{(l)}_a+{1}\right]}
\prod_{b=a+1}^{\la^{(l)}}\frac{[v^{(l)}_b-v^{(l)}_a+1]}{[v^{(l)}_b-v^{(l)}_a]}\right).\nn\\
&&\lb{entireU}
\ena

Furthermore, in order to compare with the stable envelopes, it is convenient to consider the following 
expression. See Sec.\ref{directCom}. 
\bea
\cW_I(t,z,\Pi)&=&\frac{W_I(t,z,\Pi)}{E_\la(t)}={\rm Sym}_{t^{(1)}}\cdots {\rm Sym}_{t^{(N-1)}}\ {\cU}_I(t,z,\Pi),\lb{def:cW}\\
\cU_I(t,z,\Pi)&=&\prod_{l=1}^{N-1}\frac{\prod_{a=1}^{\la^{(l)}} 
    u^{(l)}_I(t^{(l)}_a,t^{(l+1)},\Pi_{\mu_{\mbox{\tiny${i^{(l)}_a}$}},l+1}q^{-2C_{\mu_{\mbox{\tiny${i^{(l)}_a}$}},l+1}(i^{(l)}_a)})}{\prod_{1\leq a<b\leq \la^{(l)}}{[v^{(l)}_a-v^{(l)}_b]}{[v^{(l)}_b-v^{(l)}_a-1]}},\lb{def:cU}
\ena
where $t^{(l+1)}=(t^{(l+1)}_1,\cdots,t^{(l+1)}_{\la^{(l+1)}})$, and 
\bea
&&E_\la(t)=\prod_{l=1}^{N-1}\prod_{a=1}^{\la^{(l)}}\prod_{b=1}^{\la^{(l)}}[v^{(l)}_b-v^{(l)}_a+1],\lb{def:Ela}\\
&&u^{(l)}_I(t^{(l)}_a,t^{(l+1)},\Pi_{j,k})= 
\left.
\frac{\left[v^{(l+1)}_b-v^{(l)}_a+(P+h)_{j,k}
\right]}{[(P+h)_{j,k}
]}
\right|_{i^{(l+1)}_b=i^{(l)}_a}
\nn\\
&&\qquad 
\qquad\qquad\times \prod_{b=1\atop i^{(l+1)}_b>i^{(l)}_a}^{\la^{(l+1)}}{\left[v^{(l+1)}_b-v^{(l)}_a\right]}
\prod_{b=1\atop i^{(l+1)}_b<i^{(l)}_a}^{\la^{(l+1)}}{\left[v^{(l+1)}_b-v^{(l)}_a+{1}\right]}
\lb{def:u}
\ena
for $1\leq j< k\leq N$. 
 
\noindent
{\it Remark. } 
In the trigonometric $(p\to 0)$ and non-dynamical (neglecting the factors depending on $P+h$ ) limit  $W_I$ and $\cW_I$ coincide  with 
$W_I$ and $\tW_I$ discussed in \cite{RTV}, respectively. See also \cite{Mimachi,MN}.

\subsection{Properties of the elliptic weight functions}\lb{Properties}

\subsubsection{Triangular property}
For $I,J\in \cI_{\la}$, let $I^{(l)}=\{i^{(l)}_1<\cdots<i^{(l)}_{\la^{(l)}}\}$ and $J^{(l)}=\{j^{(l)}_1<\cdots<j^{(l)}_{\la^{(l)}}\}$ 
$(l=1,\cdots,N)$.
 Define a partial ordering $\leqslant$ by 
\be
I\leqslant J \Leftrightarrow i^{(l)}_a \leq j^{(l)}_a\qquad \forall l, a. 
\en
Let us denote by $t=z_I$ the speciatization $t^{(l)}_a=z_{i^{(l)}_a}$ $(l=1,\cdots,N-1, a=1,\cdots,\la^{(l)})$\cite{RTV}.
The weight function has the following triangular property.
\begin{prop}\lb{triangular}
For $I,J\in \cI_\la$, 
\begin{itemize}
\item[(1)] $\cW_{J}(z_I,z,\Pi)=0$ unless $I\leqslant J$.
\item[(2)] 
\be
&&\cW_{I}(z_I,z,\Pi)=\prod_{1\leq k<l\leq N}\prod_{a\in I_k}\left(\prod_{b\in I_l\atop a<b}{[u_b-u_a]}\prod_{b\in I_l\atop a>b}{[u_b-u_a+1]}\right).
\en
\end{itemize}

\end{prop}

For $\sigma\in \gS_n$, let us denote $\sigma^{-1}(I)=I_{\mu_{\sigma(1)}\cdots\mu_{\sigma(n)}}$ and 
$\sigma(z)=(z_{\sigma(1)},\cdots,z_{\sigma(n)})$.  
Following  \cite{RTV}, 
let us set ${\cW}_{\sigma,I}(t,z,\Pi)={\cW}_{\sigma(I)}(t,\sigma(z),\Pi)$ 
and $\cW_{\id,I}(t,z\Pi)=\cW_I(t,z,\Pi)$. 
Let us  consider the matrix $\hW_{\sigma}(z,\Pi)$, whose $(I,J)$th element is given by $\cW_{\sigma,J}(z_I,z,\Pi)$ $(I,J\in \cI_\la)$. 
We put the matrix elements in the decreasing order with respect to $\leqslant$. 
Then Proposition \ref{triangular} yields that the matrix $\hW_{\id}(z,\Pi)$ is lower triangular, whereas 
$\hW_{\sigma_0}(z,\Pi)$ with $\sigma_0$ being the longest element in $\gS_n$ is upper triangular. 
In particular, for generic $u_a\ (a=1,\cdots,n)$,  $\hW_{\sigma}(z,\Pi)$ is invertible.

\subsubsection{Transition property}\lb{TransProp}
\begin{prop}\lb{transitionProp}
Let $I=I_{\mu_1\cdots\mu_i\mu_{i+1}\cdots\mu_n}\in \cI_\la$. 
\bea
&&{\cW}
_{I_{\cdots\ \mu_{i+1}\mu_{i}\cdots}}
(t, 
\cdots,z_{i+1},z_{i},\cdots
,\Pi)\nn\\
&&=
\sum_{\mu_{i}',\mu_{i+1}'}\bar{R}(z_{i}/z_{i+1},\Pi q^{-2\sum_{j=i}^{n}<\bep_{\mu_j},h>})_{\mu_{i}\mu_{i+1}}^{\mu_{i}'\mu_{i+1}'}\ 
{\cW}_{I_{\cdots\ \mu'_i\mu'_{i+1}\cdots}}(t, 
\cdots,z_i,z_{i+1},\cdots
,\Pi).\lb{transsi}
\ena
\end{prop}
Note that since $H_\la(t,z)$ is a symmetric function in $z_1,\cdots, z_n$,  
 $\tW_I(t,z,\Pi)$ has the same property. 

\subsubsection{Orthogonality}
Noting \eqref{def:cW} and the remark in Sec.5.3 in \cite{Konno17}, where $E_\la(t,z)$ is the same as $E_\la(t)$ in \eqref{def:Ela}, we have the following property.
\begin{prop}\lb{orthogonalProp}
For $J,K\in \cI_\la$, 
\be
&&\sum_{I\in \cI_\la}\frac{\cW_J(z_I,z,\Pi^{-1}q^{2\sum_{j=1}^{n}<\bep_{\mu_j},h>})\cW_{\sigma_0(K)}(z_I,\sigma_0(z),\Pi)}{Q(z_I)R(z_I)}=\delta_{J,K},
\en
where $\sum_{j=1}^n\bep_{\mu_j}$ is the weight associated with $\la$ (Sec.\ref{combNot}), and
\be
&&Q(z_I)=\prod_{1\leq k<l\leq N}\prod_{a\in I_k}\prod_{b\in I_l}[u_b-u_a+1],\\
&&R(z_I)=\prod_{1\leq k<l\leq N}\prod_{a\in I_k}\prod_{b\in I_l}[u_b-u_a].
\en
\end{prop}
In Sec.\ref{subGeomRep}, a consistency between this property and the formula in Theorem \ref{XtW} becomes a key to obtain a  geometric representation of the elliptic quantum group.\footnote{In \cite{FRV,RTV17},  the dynamical shift  in the orthogonality relation is missing. 
}

\subsubsection{Quasi-periodicity}
Remember  that we set $t^{(l)}_a=q^{2v^{(l)}_a}$, $z_k=q^{2u_k}$ and $\Pi_{j,k}=q^{2(P+h)_{j,k}}$.  Note that  $t^{(l)}_a\mapsto pt^{(l)}_a\Leftrightarrow v^{(l)}_a\mapsto v^{(l)}_a+r$ and
$t^{(l)}_a\mapsto e^{-2\pi i}t^{(l)}_a\Leftrightarrow v^{(l)}_a\mapsto v^{(l)}_a+r\tau$. 
From \eqref{thetaquasiperiod} and Proposition \ref{combC} we obtain the following statement.   
\begin{prop}\lb{quasiperiod}
For $I\in \cI_\la$, the weight functions $\cW_I(t,z,\Pi)$ has the following quasi-periodicity. 
\be
&&\cW_I(\cdots,pt^{(l)}_a,\cdots, z,\Pi)=(-1)^{\la_{l+1}-\la_l+2}\cW_I(\cdots,t^{(l)}_a,\cdots, z,\Pi),\\
&&\cW_I(\cdots,e^{-2\pi i}t^{(l)}_a,\cdots, z,\Pi)\\
&&=(-e^{-\pi i \tau})^{\la_{l+1}-\la_l+2}\\
&&\hspace{-0.5cm}\times \exp\left\{-\frac{2\pi i}{r}\left((\la_{l+1}-\la_l)v^{(l)}_a-\sum_{b=1}^{\la^{(l+1)}}v^{(l+1)}_b
+2\sum_{b=1}^{\la^{(l)}}v^{(l)}_b
-\sum_{b=1}^{\la^{(l-1)}}v^{(l-1)}_b-(P+h)_{l,l+1}-\la_{l+1}\right)\right\}\\
&& \times \cW_I(\cdots,t^{(l)}_a,\cdots, z,\Pi) \qquad\qquad (1\leq a\leq \la^{(l)}, 1\leq l\leq N-1).
\en

\end{prop}

\vspace{2mm}
\noindent
{\it Remark.}\ 
For $\la=(\la_1,\cdots,\la_N)\in \N^N$,  let $x={}^t(x^{(1)}_1,\cdots,x^{(1)}_{\la^{(1)}},\cdots,x^{(N-1)}_1,\cdots,x^{(N-1)}_{\la^{(N-1)}})\in \C^M$, where $M=\sum_{l=1}^{N-1}\la^{(l)}=\sum_{l=1}^{N-1}(N-l)\la_l$. 
From Proposition \ref{quasiperiod} one can deduce a symmetric integral $M\times M$ matrix $N$ and a vector $\xi\in (\C/r\Z)^M$, which imply  
 the following quadratic form $N(x)={}^tx N x$ and the linear form $\xi(x)={}^tx \xi$. 
\be
&&N(x)=-\sum_{l=1}^{N-1}\sum_{a=1}^{\la^{(l)}}\sum_{b=1}^{\la^{(l)}}\left(x^{(l)}_a-x^{(l)}_b\right)^2+
\sum_{l=1}^{N-2}\sum_{a=1}^{\la^{(l)}}\sum_{b=1}^{\la^{(l+1)}}\left(x^{(l)}_a-x^{(l+1)}_b\right)^2
+n\sum_{a=1}^{\la^{(N-1)}}(x^{(N-1)}_a)^2,\\
&&\xi(x)=-\sum_{l=1}^{N-1}\sum_{a=1}^{\la^{(l)}}x^{(l)}_a\left((P+h)_{l,l+1}+\la_{l+1}\right)-\sum_{a=1}^{\la^{(N-1)}}\sum_{k=1}^nx^{(N-1)}_au_k.
\en
Then by Appel-Humbert theorem\cite{LB}, a pair $(N,\xi)$  characterizes a line bundle $\cL(N,\xi)\ :(\C^M\times \C)/\Lambda^M\to \C^M$, where 
$\Lambda=r\Z+r\Z \tau$, with action 
\be
\omega\cdot(x,\eta)=(x+\omega,e_{\omega}(x)\eta),\qquad \omega\in \Lambda^M,\ x\in \C^M,\ \eta\in \C,
\en
and cocycle 
\be
e_{nr+mr\tau}(x)=(-1)^{{}^tnNn}(-e^{i\pi \tau})^{{}^tmNm}e^{\frac{2\pi i}{r}{}^tm(Nx+\xi)},\qquad n,m\in \Z^M.
\en 
Hence ${\rm Span}_\C\{\ \cW_I(t,z,\Pi)\ (I\in \cI_\la)\ \}$ is a space of meromorphic sections of $\cL(N,\xi)$.

\subsubsection{Shuffle algebra structure}

For $\la=(\la_1,\cdots,\la_N)\in \N^N$, $|\la|=n$, let $z^{(n)}=(z_1,\cdots,z_n)\in (\C^*)^n$.   
For $I=I_{\mu_1\cdots \mu_n}\in \cI_{\la}$, we denote by $\Pi_I$ a set of dynamical parameters $\{\Pi_{\mu_k,j}=q^{2(P+h)_{\mu_k,j}}\ (k=1,\cdots, n, j=\mu_k+1,\cdots,N) \}$,where  ${(P+h)_{j,k}}\in \C/r\Z\ (1\leq j<k\leq N)$, and set $\Pi_\la=\cup_{I\in\cI_\la}\Pi_I$. 
\begin{dfn}
For $\la=(\la_1,\cdots,\la_N)\in \N^N$, $|\la|=n$, 
we define $\cM^{(n)}_\la(z^{(n)},\Pi_\la)$ to be  the space of moromorphic functions $F(t;z,\Pi)$ of $M$ variables  $t=(t^{(1)}_1,\cdots,t^{(1)}_{\la^{(1)}},\cdots,$  $ t^{(N-1)}_{1},\cdots,t^{(N-1)}_{\la^{(N-1)}})$ such that 
\begin{itemize}
\item[(1)] $F(t;z,\Pi)$ is symmetric in $t^{(l)}_{1},\cdots,t^{(l)}_{\la^{(l)}}$ for each $l\in \{1,\cdots,N-1\}$. 
\item[(2)] ${F}(t;z,\Pi)$ has the quasi-periodicity
\be
&&{F}(\cdots, pt^{(l)}_a,\cdots;z,\Pi )={F}(t;z,\Pi),\\
&&{F}(\cdots, e^{-2\pi i}t^{(l)}_a;\cdots,z,\Pi )=\exp\left\{\frac{2\pi i}{r}\left((P+h)_{l,l+1}-\la_{l}\right)\right\}{F}(t;z,\Pi)
\en
$(l=1,\cdots,N-1,\ a=1,\cdots,\la^{(l)})$.

\end{itemize}
\end{dfn}

Let us consider the subspace space $\cM^{+(n)}_\la(z^{(n)},\Pi_\la):={\rm Span}_\C\{\ \widetilde{W}_I(t,z,\Pi)\ (I\in \cI_\la)\ \}$ 
of $\cM^{(n)}_\la(z,\Pi_\la)$. 
From Proposition \ref{triangular}, we obtain 
\begin{prop} 
 $\ds{{\rm dim}_\C \cM^{+(n)}_\la(z^{(n)},\Pi_\la)=\frac{n!}{\la_1!\cdots \la_N!}}$. 
\end{prop}

Consider a graded $\C$-vector space
\be
&&\cM(z,\Pi)=
\bigoplus_{n\in\N}\bigoplus_{\la\in \N^N\atop |\la|=n}\cM^{(n)}_\la(z^{(n)},\Pi_\la)
\en
with $\cM^{(0)}_{(0,\cdots,0)}(z^{(0)},\Pi)=\C 1$.  

\begin{dfn}\lb{defstarprod}
For $F(t;z^{(m)}, \Pi_I)\in \cM^{(m)}_\la(z^{(m)},\Pi_\la)$,  $G(t';{z'}^{(n)},\Pi'_{I'})\in {\cM^{(n)}}_{\la'}({z'}^{(n)},\Pi'_{\la'})$, 
we define the bilinear product $\star$ on $\cM(z,\Pi)$ by 
\bea
&&(F\star G)(t^{(1)}_1,\cdots,  t^{(1)}_{\la^{(1)}+{\la'}^{(1)}}, \cdots, t^{(N-1)}_1,\cdots,  t^{(N-1)}_{\la^{(N-1)}+{\la'}^{(N-1)}};  z_1,\cdots,z_{m+n},
{\Pi}_{I+I'}
)\nn\\
&&:=\frac{1}{\prod_{l=1}^{N-1}\la^{(l)}!\la^{'(l)}!}{\rm Sym}^{(1)}\cdots {\rm Sym}^{(N-1)} \left[
F(t,z,\Pi_I q^{-2\sum_{j=1}^n<\bep_{\mu_j'},h>})\ G(t',z',\Pi'_{I'})\ {\Xi}(t,t',z,z')
\right], \nn\\
&&\lb{shufflprod}
\ena
where $I'=I'_{\mu_1'\cdots\mu'_n}$ and 
\be
&&{\Xi}(t,t',z,z')=\prod_{l=1}^{N-1}
\prod_{a=1}^{\la^{(l)}}\left(\prod_{b=1}^{\la^{'(l+1)}}\frac{[{v'_b}^{(l+1)}-v^{(l)}_a]}{[{v'_b}^{(l+1)}-v^{(l)}_a+1]}
\prod_{c=1}^{\la^{'(l)}}\frac{[{v'_c}^{(l)}-v^{(l)}_a+1]}{[{v'_c}^{(l)}-v^{(l)}_a]}\right). 
\en
In the LHS of \eqref{shufflprod},  we set $t^{(l)}_{\la^{(1)}+a}:={t'}^{(l)}_a\ (a=1,\cdots, {\la'}^{(l)})$, $z_{m+k}:=z'_{k}\ (k=1,\cdots,n)$ and 
 $\Pi_{I+I'}=\{\Pi_{\mu_k,j}\ (k=1,\cdots,m+n, j=\mu_k+1,\cdots,N) \}$, where 
 ${\Pi}_{\mu_{m+k},j}:=\Pi'_{\mu'_k,j}$ $(k=1,\cdots,n, j=\mu'_k+1,\cdots,N )$. 
\end{dfn}
This endows $\cM(z,\Pi)$ with a structure of an associative unital algebra with the unit $1$. 
In \cite{FRV}, a $\slt$ version of the $\star$-product is given. 

Let us consider the subspace of $\cM(z,\Pi)$. 
\be
&&\cM^+(z,\Pi)=\bigoplus_{n\in \N}\bigoplus_{\la\in \N^N\atop |\la|=n } \cM^{+(n)}_\la(z^{(n)},\Pi_\la). 
\en 
All the elements in ${\cM}^+(z,\Pi)$ satisfy the following pole and wheel conditions. 
For $F(t;z,\Pi)\in \cM^{+(n)}_\la(z^{(n)},\Pi_\la)$, 
\begin{itemize}
\item[1)] there exists an entire function $f(t;z,\Pi)\in \Theta^+_\la(z^{(n)},\Pi_\la)={\rm Span}_\C\{\ W_I(t,z^{(n)},\Pi)\ (I\in \cI_\la)\ \}
$ such that 
\be
&& F(t;z,\Pi)=\frac{f(t;z,\Pi)}{H_\la(t,z)}.
\en
\item[2)] $f(t;z,\Pi)=0$ once $t^{(l)}_a/t^{(l+\vep)}_c=q^{2\vep}$ and $t^{(l+\vep)}_c/t^{(l)}_b=1$ for some 
$l,\vep, a,b,c$, where $\vep\in\{\pm 1\}$, $l=1,\cdots,N$, $a,b=1,\cdots,\la^{(l)}$, $c=1,\cdots,\la^{(l+\vep)}$ and $t^{(N)}_a=z_a$.

\end{itemize}

\begin{prop}\lb{shuffleMp}
The subspace $\cM^+(z,\Pi)\subset \cM(z,\Pi)$ is $\star$-closed.  
\end{prop}

\section{Finite Dimensional Representations
}\lb{finRep}
In this section we construct finite dimensional tensor product  representations of the elliptic quantum group 
$E_{q,p}(\glnh)$ and  $U_{q,p}(\glnh)$ on  the  
Gelfand-Tsetlin basis. 

\subsection{Finite dimensional tensor product representations
}
Let $(\pi_z,\hV_{z})$ denote the $N$-dimensional dynamical evaluation representation of $E_{q,p}(\glnh)$:
 $\hV_z=\hV[z,z^{-1}]$ with $\hV=\oplus_{\mu=1}^N\FF v_{\mu}$. The level-0 action of the  $L$-operator  $\hL^\pm(z)
$ or the dynamical $L$-operator $L^\pm(z,P)$ introduced in Sec.\ref{dynL}
is given by
\be
\pi_{z}(\hL^\pm_{ij}(1/w))v_\nu
=\pi_{z}(L^\pm_{ij}(1/w,P)e^{-Q_{\bep_j}})v_\nu
=\sum_{\mu=1}^N\bar{R}(z/w,\Pi^*)_{i\mu}^{j\nu}v_\mu
\en
with $e^{Q_\al}v_\mu=v_{\mu}\ (\al\in \bh^*)$, where $\Pi^*_{j,l}=q^{2P_{j,l}}$ as before.

The action on the tensor product space is obtained by the co-algebra structure presented in Appendix \ref{HHopfalgebroid}. 
\begin{prop}\lb{coproL}
$\hL^\pm(1/w)$ acts on $\hV_w \tot  \hV_{z_1} \tot \cdots \tot \hV_{z_n}$ by 
\be
&&(\pi_{z_1}\otimes \cdots\otimes \pi_{z_n}){\Delta'}^{(n-1)}(\hL^\pm(1/w)) \\
&&=\bar{R}^{(0n)}(z_n/w,\Pi^* q^{2\sum_{j=1}^{n-1}h^{(j)}}){\bR}^{(0n-1)}(z_{n-1}/w,\Pi^* q^{2\sum_{j=1}^{n-2}h^{(j)}})\cdots \bar{R}^{(01)}(z_1/w,\Pi^*).
\en
\end{prop}
\noindent
{\it Proof. } It is enough to show  the $n=2$ case. 
\be
(\pi_{z_1}\otimes \pi_{z_2}){\Delta'}(\hL^\pm_{ij}(1/w)) v_{\mu}\tot v_{\nu}
&=&\sum_{k}\pi_{z_1}(\hL^\pm_{kj}(1/w))v_\mu\tot \pi_{z_2}(\hL^\pm_{ik}(1/w))v_\nu\\
&=&\sum_{k,\mu',\nu'}\bR(z_1/w,\Pi^*)_{k\mu'}^{j\mu}v_{\mu'}\tot \bR(z_2/w,\Pi^*)_{i\nu'}^{k\nu}v_{\nu'}\\
&=&\sum_{\mu',\nu'}\sum_{k}\bR(z_2/w,\Pi^* q^{2h^{(1)}})_{i\nu'}^{k\nu}\bR(z_1/w,\Pi^*)_{k\mu'}^{j\mu}v_{\mu'}\tot v_{\nu'}\\
&=&\sum_{\mu',\nu'}\left(\bR^{(02)}(z_2/w,\Pi^* q^{2h^{(1)}})\bR^{(01)}(z_1/w,\Pi^*)\right)_{i\mu'\nu'}^{j\mu\nu}v_{\mu'}\tot v_{\nu'}.
\en
To obtain the third equality we used \eqref{AtotB}.
\qed

It is also useful to write down the comultiplication fomula of the dynamical $L$-operator, which is equivalent to Proposition \ref{coproL}.
\begin{prop}\lb{coproDL}
The dynamical $L$-operator  $L^+(1/w,P)$ acts on  $\hV_w\otimes \hV_{z_1}\otimes \cdots \otimes \hV_{z_n}$, where 
 $\otimes$ denotes the usual tensor product, by 
 \be
&&(\pi_{z_1}\otimes \cdots\otimes \pi_{z_n}){\Delta'}^{(n-1)}(L^+_{ij}(1/w,P)) \\
&&= \sum_{k_1,\cdots,k_{n-1}=1}^NL^+_{k_1j}(z_1/w,P)\otimes L^+_{k_2k_1}(z_1/w,P+h^{(1)})\otimes \cdots \otimes 
L^+_{ik_{n-1}}(z_n/w,P+\sum_{j=1}^{n-1}h^{(j)}).
\en
\end{prop}

\subsection{The Gelfand-Tsetlin basis}

\begin{dfn}
The Gelfand-Tsetlin basis is a basis of the level-$0$ representation of $E_{q,p}(\glnh)$ or $U_{q,p}(\glnh)$ 
consisting of the simultaneous  eigenvectors of the Gelfand-Tsetlin subalgebra $\mathfrak{G}$ in Definition \ref{def:GTsubalg}.  
\end{dfn}

We consider the  the Gelfand-Tsetlin (GT) basis in $\hV_{z_1}\tot \cdots \tot \hV_{z_n}$.
 Following \cite{RTV} we construct it as follows.
Firstly we realize $\goth{S}_n$  in terms of the elliptic dynamical $R$ matrix in \eqref{ellR}. 
Let define $\widetilde{S}_i(P)$ by 
\be
&&\widetilde{S}_i(P):=\cP^{(i i+1)}\bR^{(i i+1)}(z_{i}/z_{i+1},\Pi^* q^{2{\sum_{j=1}^{i-1}h^{(j)}}})s^z_i, 
\en
where
\be
&&\cP: v\tot w \mapsto w\tot v, \qquad s^z_i f(\cdots,z_i ,z_{i+1},\cdots )=f(\cdots,z_{i+1}, z_i ,\cdots )
\en
Then by using the dynamical Yang-Baxter equation \eqref{DYBE} and the unitarity relation \eqref{unitarity} one can show the following.
\begin{prop}
\be
&&\hspace{-2cm}\tS_i(P)\tS_{i+1}(P)\tS_i(P)=\tS_{i+1}(P)\tS_{i}(P)\tS_{i+1}(P),\\
&&\hspace{-2cm}\tS_i(P)\tS_j(P)= \tS_j(P)\tS_i(P) \qquad\qquad (|i-j|>1)\\
&&\hspace{-2cm}\tS_{i}(P)^2=1.
\en
\end{prop}

For $\la\in \N^N, |\la|=n$, $I=I_{\mu_1\cdots \mu_n}\in \cI_{\la}$, we set 
\be
&&v_I=v_{\mu_1\cdots\mu_n}:=v_{\mu_1}\tot \cdots \tot v_{\mu_n}.
\en 
We define the Gelfand-Tsetlin basis $\{\xi_I\}_{I\in \cI_{\la}}$ by
\bea
&&\hspace{-2cm} \xi_{I^{max}}:=v_{I^{max}},\qquad 
\xi_{s_i(I)}:=\tS_i(P)\xi_I,\lb{def:GT}
\ena  
where 
\be
&& I^{max}=I_{\tiny\underbrace{N\cdots N}_{\la_N}\ \cdots \underbrace{1\cdots 1}_{\la_1}}. 
\en

Let us consider the change of basis matrix $\widehat{X}=(X_{IJ}(z,P))_{I,J\in \cI_{\la}}$:
\bea
&&\xi_I=\sum_{J\in \cI_{\la}}X_{IJ}(z,P)v_J.\lb{xiv}
\ena
 Here we put the matrix elements  in the decreasing order $I^{max}\geqslant\cdots \geqslant I^{min}$. 
Then by construction, 
$\widehat{X}$ is a lower triangular matrix.  Furthermore the following remarkable relationship between $\widehat{X}$ 
and the specialized elliptic weight functions becomes a key to obtain a geometric interpretation of the results 
in the next subsection. See Sec.\ref{subGeomRep}.   
\begin{thm}\lb{XtW}
\bea
&&X_{IJ}(z,P)=\tW_J(z^{-1}_I,z^{-1},\Pi^* q^{2\sum_{j=1}^n<\bep_{\mu_j},h>}).\lb{XW}
\ena
\end{thm}
\noindent
{\it Proof.}\ Let $J=I_{\mu_1\cdots\mu_i\mu_{i+1}\cdots\mu_n}\in \cI_{\la}$. By definition,
\be
\xi_{s_i(I)}&=&\sum_{J}X_{s_i(I)J}(z,P)v_J\\
&=&\tS_i(P)\xi_I=\sum_{J}X_{IJ}(s_i(z),P)\tS_i(P)v_J\\
&=&\sum_{J, \mu'_i, \mu'_{i+1}}X_{IJ}(s_i(z),P)\bR(z_i/z_{i+1},\Pi^* q^{2\sum_{j=1}^{i-1}<\bep_{\mu_j},h>})_{\mu'_i \mu'_{i+1}}^{\mu_i\mu_{i+1}}v_{\mu_1}\tot \cdots\tot v_{\mu'_i}\tot v_{\mu'_{i+1}}\tot \cdots \tot v_{\mu_n}
.
\en
Hence we obtain 
\bea
&&X_{s_i(I)J}(z,P)=X_{IJ}(s_i(z),P)\lb{rec1}
\ena
for $\mu_i=\mu_{i+1}$, and
\bea
&&\left(X_{s_i(I)J}(z,P)\ X_{s_i(I)s_i(J)}(z,P)\right)\nn\\
&&=\left(X_{IJ}(s_i(z),P)\ X_{Is_i(J)}(s_i(z),P)\right)\cP_2{}^t\bR(z_i/z_{i+1},\Pi^* q^{2\sum_{j=1}^{i-1}<\bep_{\mu_j},h>})_{\mu_i,\mu_{i+1}}
\lb{rec2}
\ena
for $\mu_i>\mu_{i+1}$. Here we set
\bea
&&\cP_2=\mat{0&1\cr 1&0\cr},\qquad \bR(z,\Pi^* )_{\mu_i,\mu_{i+1}}
=\left(\mmatrix{\bR(z,\Pi^* )_{\mu_{i+1}\mu_i}^{\mu_{i+1}\mu_i}& \bR(z,\Pi^* )_{\mu_{i+1}\mu_i}^{\mu_i\mu_{i+1}}\cr
\bR(z,\Pi^* )_{\mu_{i}\mu_{i+1}}^{\mu_{i+1}\mu_i}& \bR(z,\Pi^* )_{\mu_{i}\mu_{i+1}}^{\mu_{i}\mu_{i+1}}}\right).\lb{bRmat}
\ena 
Note that \eqref{rec1} and \eqref{rec2} determine the whole matrix elements in $\widehat{X}$ recursively  starting from 
$X_{I^{max}I^{max}}(z,P)=1$. 

On the other hand, from Proposition \ref{transitionProp} with replacing $\Pi$ by $\Pi^*$ we have 
\bea
&&\tW_J(t,s_i(z),\Pi^*)=\tW_J(t,z,\Pi^*)\lb{relW1}
\ena
if $\mu_i=\mu_{i+1}$, and 
\be
&&\left(\tW_{J}(t,s_i(z),\Pi^*)\ \tW_{s_i(J)}(t,s_i(z),\Pi^*)\right)\nn\\
&&=\left(\tW_{J}(t,z,\Pi^*)\ \tW_{s_i(J)}(t,z,\Pi^*)\right)\cP_2{}^t\bR(z_i/z_{i+1},\Pi^* q^{-2\sum_{j=i}^{n}<\bep_{\mu_j},h>})_{\mu_i,\mu_{i+1}}
\en
if $\mu_i\not=\mu_{i+1}$. 
Using 
\be
&&\left(\cP_2{}^t\bR(z,\Pi^*)_{\mu_i,\mu_{i+1}}\right)^{-1}=\cP_2{}^t\bR(z^{-1},\Pi^*)_{\mu_i,\mu_{i+1}},
\en
we obtain in particular  for $\mu_i>\mu_{i+1}$
\bea
&&\left(\tW_{J}(t,z,\Pi^*)\ \tW_{s_i(J)}(t,z,\Pi^*)\right)\nn\\
&&=\left(\tW_{J}(t,s_i(z),\Pi^*)\ \tW_{s_i(J)}(t,s_i(z),\Pi^*)\right)\cP_2{}^t\bR(\left(z_i/z_{i+1}\right)^{-1},\Pi^* q^{-2\sum_{j=i}^{n}<\bep_{\mu_j},h>})_{\mu_i,\mu_{i+1}}
\lb{relW2}
\ena
Specializing $t=s_i(z)_I$ and noting 
\be
&&\tW_{J}(s_i(z)_I,z,\Pi^*)=\tW_{J}(z_{s_i(I)},z,\Pi^*) 
\en
etc.,  we obtain from \eqref{relW1} and \eqref{relW2}
\bea
&&\tW_J(s_i(z)_I,s_i(z),\Pi^*)=\tW_J(z_{s_i(I)},z,\Pi^*)\lb{srelW1}
\ena
if $\mu_i=\mu_{i+1}$, and 
\bea
&&\left(\tW_{J}(z_{s_i(I)},z,\Pi^* )\ \tW_{s_i(J)}(z_{s_i(I)},z,\Pi^* )\right)\nn\\
&&=\left(\tW_{J}(s_i(z)_I,s_i(z),\Pi^*)\ \tW_{s_i(J)}(s_i(z)_I,s_i(z),\Pi^*)\right)\cP_2{}^t\bR(\left(z_i/z_{i+1}\right)^{-1},\Pi^* q^{-2\sum_{j=i}^{n}<\bep_{\mu_j},h>})_{\mu_i,\mu_{i+1}}\nn\\
&&\lb{srelW2}
\ena
if $\mu_i>\mu_{i+1}$. 
Therefore one finds that $\tW_J(z^{-1}_I,z^{-1},\Pi^* q^{2\sum_{j=1}^n<\bep_{\mu_j},h>})$  satisfy the same recursion relations  
as  \eqref{rec1} and \eqref{rec2} for $X_{IJ}(z,P)$. In addition their initial conditions coincide: $\tW_{I^{max}}(z^{-1}_{I^{max}},z^{-1},\Pi^* q^{2\sum_{j=1}^n<\bep_{\mu_j},h>})=1=X_{I^{max}I^{max}}(z,P)$. \qed

\vspace{2mm}

\noindent
{\it Example. } The case $N=2, n=3, \la=(2,1)$. We have $\cI_\la=\{ I_{211}\geqslant I_{121}\geqslant I_{112}\}$ and 
\be
\mat{\xi_{211}\cr 
\xi_{{121}}\cr 
\xi_{112}\cr} =\mat{1&0&0\cr
c(u_{1,2},P_{1,2})&\bb(u_{1,2})&0\cr
c(u_{1,3},P_{1,2})&\bb(u_{1,3})c(u_{2,3},P_{1,2}+1)&\bb(u_{1,3})\bb(u_{2,3})\cr}
\mat{v_{211}\cr 
v_{{121}}\cr 
v_{112}\cr},
\en  
where $u_{i,j}=u_i-u_j$. 
On the other hand we have
\be
\hW_{\id}(z,\Pi^*)&=&
\mat{\tW_{I_{211}}(z_{I_{211}},z,\Pi^*)&0&0\cr
\tW_{I_{211}}(z_{I_{121}},z,\Pi^*)&
\tW_{I_{121}}(z_{I_{121}},z,\Pi^*)&0\cr
\tW_{I_{211}}(z_{I_{112}},z,\Pi^*)&
\tW_{I_{121}}(z_{I_{112}},z,\Pi^*)&\tW_{I_{112}}(z_{I_{112}},z,\Pi^*)\cr}\\
&=&\mat{1&0&0\cr
\frac{[P_{1,2}-1+u_{2,1}][1]}{[u_{2,1}+1][P_{1,2}-1]}&
\frac{[u_{2,1}]}{[u_{2,1}+1]}&0\cr
\frac{[P_{1,2}-1+u_{3,1}][1]}{[u_{3,1}+1][P_{1,2}-1]}&
\frac{[u_{3,1}]}{[u_{3,1}+1]}\frac{[P_{1,2}+u_{3,2}][1]}{[u_{3,2}+1][P_{1,2}]}&
\frac{[u_{3,1}]}{[u_{3,1}+1]}\frac{[u_{3,2}]}{[u_{3,2}+1]}\cr}.
\en

\subsection{Action of the elliptic currents}
In order to derive an action of the elliptic quantum group on the GT basis $\{\xi_I\}$, 
the following property of the symmetrization operators $\widetilde{S}_i(P)$ is useful \cite{RTV}.   
\begin{prop}
\be
&&\tS_i(P){\Delta'}^{(n-1)}(\hL^\pm(w))={\Delta'}^{(n-1)}(\hL^\pm(w))\tS_i(P+h^{(0)})
\en 
\end{prop}
\noindent
{\it Proof.} Use the dynamical Yang-Baxter equation. 
\qed

\noindent
Thanks  to this proposition it suffices to construct an action of 
${\Delta'}^{(n-1)}(\hL^\pm(w))$ on $\xi_{I^{max}}$.

From Theorem \ref{HCqminor} and Proposition \ref{coproL} we obtain the following level-0 action of the half-currents of $E_{q,p}(\glnh)$ 
on the GT basis. Note that  at the level 0 we have $L^-(w,P)=L^+({p}w,P)$, hence   
$\cK^-_j(w)=\cK^+_j(pw), E^-_{j+1,j}(w,P)=E^+_{j+1,j}(pw,P),\ F^-_{j,j+1}(w,P)=F^+_{j,j+1}(pw,P)$.
\begin{thm}\lb{actHC}
Let $\{\xi_I\ |\ I\in \cI_\la, \la=(\la_1,\cdots,\la_N)\in \N^N, |\la|=n\}$ be the GT basis of $\hV_{z_1}\tot \cdots \tot \hV_{z_n}$. 
Under the abbreviation $\cK^\pm_j(1/w)=(\pi_{z_1}\otimes \cdots \otimes \pi_{z_n}){\Delta'}^{(n-1)}(\cK^\pm_j(1/w))$, $E^\pm_{j+1,j}(1/w,P)=$\\ 
$(\pi_{z_1}\otimes \cdots \otimes \pi_{z_n}){\Delta'}^{(n-1)}(E^\pm_{j+1,j}(1/w,P))$ and  
$F^\pm_{j,j+1}(1/w,P)=(\pi_{z_1}\otimes \cdots \otimes \pi_{z_n}){\Delta'}^{(n-1)}(F^\pm_{j,j+1}(1/w,P))$, we have
\bea
&&\cK_j^\pm(1/w)\xi_I=\left.\prod_{k=1}^{j-1}\prod_{a\in I_k}\frac{[u_a-v]}{[u_a-v+1]}\right|_{\pm}\left.\prod_{l=j+1}^N\prod_{b\in {I_l}}\frac{[u_b-v-1]}{[u_b-v]}\right|_{\pm}\ \xi_I,\lb{actK}\\
&&\hspace{-0.2cm}E^\pm_{j+1,j}(1/w,P)\xi_I=\sum_{i\in I_{j+1}}\left.\frac{[P_{j,j+1}-u_i+v][1]}{[P_{j,j+1}][u_i-v]}\right|_{\pm}\prod_{k\in I_{j+1}\atop \not=i}\frac{[u_i-u_k+1]}{[u_i-u_k]}\ \xi_{I^{i'}},\lb{actE}\\
&&\hspace{-0.7cm}F^\pm_{j,j+1}(1/w,P)\xi_I=\sum_{i\in I_{j}}\left.\frac{[P_{j,j+1}+\la_j-\la_{j+1}+u_i-v-1][1]}{[P_{j,j+1}+\la_j-\la_{j+1}-1][u_i-v]}\right|_{\pm}\prod_{k\in I_{j}\atop \not=i}\frac{[u_k-u_i+1]}{[u_k-u_i]}\ \xi_{I^{'i}},\lb{actF}
\ena
where $w=q^{2v}, z_i=q^{2u_i} \ (i=1,\cdots, n)$,  $I=(I_1,\cdots, I_N)$, and $I^{i'}\in \cI_{(\la_1,\cdots,\la_j+1,\la_{j+1}-1,\cdots,\la_N)}$ and $I^{'i}\in \cI_{(\la_1,\cdots,\la_j-1,\la_{j+1}+1,\cdots,\la_N)}$ are defined by  
\be
&&(I^{i'})_j=I_j\cup \{i\},\quad (I^{i'})_{j+1}=I_{j+1}-\{i\},\quad (I^{i'})_k=I_k \ \ (k\not=j,j+1), \\
&&(I^{'i})_j=I_j- \{i\},\quad (I^{'i})_{j+1}=I_{j+1}\cup\{i\},\quad (I^{'i})_k=I_k \ \ (k\not=j,j+1).
\en
The symbols $|_{\pm}$ distinguish contributions from the plus- and the minus-half-currents and specify their expansion directions as 
\be
\left.\frac{[s+v]}{[s][v]}\right|_{+}&=&w^{\frac{s}{r}}\frac{\Theta_p(q^{2s}w)}{\Theta_p(q^{2s})\Theta_p(w)} \qquad\qquad \mbox{expand in $w$}\\
\left.\frac{[s+v]}{[s][v]}\right|_{-}&=&({p}w)^{\frac{s}{r}}\frac{\Theta_p({p}q^{2s}w)}{\Theta_p(q^{2s})\Theta_p({p}w)} \\
&=&- w^{\frac{s}{r}}\frac{\Theta_p(q^{-2s}/w)}{\Theta_p(q^{-2s})\Theta_p(1/w)}\qquad \mbox{expand in $1/w$.}
\en

\end{thm}
A proof is given in Appendix \ref{ProofacHC}. 

\noindent
{\it Example.} Let us consider the case $N=3$, $n=5$, $\la=(2,2,1)$. We use the abbreviation  $\xi_{\mu_1\mu_2\mu_3\mu_4\mu_5}=\xi_{I_{\mu_1\mu_2\mu_3\mu_4\mu_5}}$, 
 as well as $\hL^+_{ij}(1/w)=(\pi_{z_1}\otimes \cdots \otimes \pi_{z_n}){\Delta'}^{(4)}(\hL^+_{ij}(1/w))$ on $\hV_{z_1}\tot \cdots\tot \hV_{z_4}$. 
Noting $I^{max}=I_{32211}\in \cI_\la$, let us consider the action on $\xi_{32211}$. 
\begin{itemize}
\item $\cK^+_3(1/w)$: from Theorem \ref{HCqminor} and Proposition \ref{coproL}, we obtain 
\be
&&\cK^+_3(1/w)\xi_{32211}=\hL^+_{33}(1/w)\xi_{32211}=\bb(u_2-v)\bb(u_3-v)\bb(u_4-v)\bb(u_5-v)\xi_{32211}.
\en

\item $E^+_{3,2}(1/w,P)$: similarly we have
\be
&&\hL^+_{32}(1/w)\xi_{32211}=\bac(u_1-v,P_{2,3})\bb(u_2-v)\bb(u_3-v)\bb(u_4-v)\bb(u_5-v)\xi_{22211}.
\en
Note that $I_{22211}$ is the maximal partition in $\cI_{(2,3,0)}$. We also obtain
\be
&&\hL^+_{33}(1/w)\xi_{22211}=\bb(u_1-v)\bb(u_2-v)\bb(u_3-v)\bb(u_4-v)\bb(u_5-v)\xi_{22211}.
\en
Hence
\be
E^+_{3,2}(1/w,P)\xi_{32211}&=&\hL^+_{33}(1/w)^{-1}\hL^+_{32}(1/w)\xi_{32211}=\frac{c(u_1-v,P_{2,3})}{\bb(u_1-v)}\xi_{22211}.
\en

\item $F^+_{2,3}(1/w,P)$: similaly we obtain 
\be
\hL^+_{23}(1/w)\xi_{32211}&=&c(u_2-v,P_{2,3}-1)\bb(u_4-v)\bb(u_5-v)\xi_{33211}\\
&&+\bb(u_2-v)c(u_3-v,P_{2,3})\bb(u_4-v)\bb(u_5-v)v_{32311}.
\en
Again note that $I_{33211}$ is the maximal partition in $\cI_{(2,1,2)}$. 
From \eqref{def:GT} we have
\be
&&v_{32311}=\frac{1}{\bb(u_{23})}\xi_{32311}-\frac{c(u_{23},P_{2,3}-1)}{\bb(u_{23})}\xi_{33211}.
\en
Substituting this and using  the identity
\be
&&c(u_2-v,P_{2,3}-1)-\bb(u_2-v)\frac{c(u_{23},P_{2,3}-1)}{\bb(u_{23})}c(u_3-v,P_{2,3})=\frac{c(u_2-v,P_{2,3})}{\bb(u_{32})}\bb(u_3-v),
\en
we obtain 
\be
\hL^+_{23}(1/w)\xi_{32211}&=&\frac{c(u_2-v,P_{2,3})}{\bb(u_{32})}\bb(u_3-v)\bb(u_4-v)\bb(u_5-v)
\xi_{33211}\\
&&+\bb(u_2-v)\frac{c(u_3-v,P_{2,3})}{\bb(u_{23})}\bb(u_4-v)\bb(u_5-v)
\xi_{32311},
\en
where $u_{ij}=u_i-u_j$. Hence 
\be
F^+_{2,3}(1/w,P)\xi_{32211}&=&\hL^+_{23}(1/w)\hL^+_{33}(1/w)^{-1}\xi_{32211}\\
&=&\frac{c(u_2-v,P_{2,3})}{\bb(u_2-v)}\frac{1}{\bb(u_{32})}
\xi_{33211}
+\frac{c(u_3-v,P_{2,3})}{\bb(u_3-v)}\frac{1}{\bb(u_{23})}
\xi_{32311}.
\en

\item $\cK^+_2(1/w)$: let us consider the action
\be
\hL^+_{22}(1/w)\xi_{32211}&=&b(u_1-v,P_{2,3})\bb(u_4-v)\bb(u_5-v)\xi_{32211}\\
&&+\bac(u_1-v,P_{2,3})c(u_2-v,P_{2,3}+1)\bb(u_4-v)\bb(u_5-v)v_{23211}\\
&&+\bac(u_1-v,P_{2,3})\bb(u_2-v)c(u_3-v,P_{2,3}+2)\bb(u_4-v)\bb(u_5-v)v_{22311}. 
\en
On the other hand  we have
\be
&&F^+_{2,3}(1/w,P)\cK^+_3(1/w)E^+_{3,2}(1/w,P)\xi_{32211}\\
&&=c(u_1-v,P_{2,3})\frac{\bac(u_1-v,P_{2,3})}{\bb(u_1-v)}\bb(u_4-v)\bb(u_5-v)\xi_{32211}\\
&&+\bac(u_1-v,P_{2,3})c(u_2-v,P_{2,3}+1)\bb(u_4-v)\bb(u_5-v)v_{23211}\\
&&+\bac(u_1-v,P_{2,3})\bb(u_2-v)c(u_3-v,P_{2,3}+2)\bb(u_4-v)\bb(u_5-v)v_{22311}. 
\en
Therefore 
\be
\cK^+_2(1/w)\xi_{32211}&=&\left(\hL^+_{22}(1/w)-F^+_{2,3}(1/w,P)\cK^+_3(1/w)E^+_{3,2}(1/w,P)\right)\xi_{32211}\\
&=&\frac{\bb(u_4-v)\bb(u_5-v)}{\bb(-u_1+v)}\xi_{32211}.
\en
The last equality follows from the identity
\be
&&b(u_1-v,P_{2,3})-c(u_1-v,P_{2,3})\frac{\bac(u_1-v,P_{2,3})}{\bb(u_1-v)}=\frac{1}{\bb(-u_1+v)}.
\en
\end{itemize}
\qed

Furthermore  noting  the formula
\bea
&&\left.\frac{[s+u]}{[s][u]}\right|_{+}-\left.\frac{[s+u]}{[s][u]}\right|_{-}=\frac{1}{[0]'}\delta(w) \lb{pmdelta}
\ena
and using  Propositions \ref{EcHc} and \ref{dynamicalHC} we obtain the following level-0 
action of the elliptic currents of $U_{q,p}(\slnh)$ on the GT basis. 
\begin{cor}\lb{actUqp}
\be
&&H_j^\pm(q^{j-N+1}/w)\xi_I=\varrho
\left.\prod_{a\in I_j}\frac{[u_a-v+1]}{[u_a-v]}\right|_{\pm}\left.\prod_{b\in {I_{j+1}}}\frac{[u_b-v-1]}{[u_b-v]}\right|_{\pm}\ e^{-Q_{\al_j}}\xi_I,\\
&&E_{j}(q^{j-N+1}/w)\xi_I=\frac{\mu^*[1]}{[0]'}
\sum_{i\in I_{j+1}}\delta(z_i/w)\prod_{k\in I_{j+1}\atop \not=i}\frac{[u_i-u_k+1]}{[u_i-u_k]}\ e^{-Q_{\al_j}}\xi_{I^{i'}},\\
&&F_{j}(q^{j-N+1}/w)\xi_I=\frac{\mu[1]}{[0]'}
\sum_{i\in I_{j}}\delta(z_i/w)\prod_{k\in I_{j}\atop \not=i}\frac{[u_k-u_i+1]}{[u_k-u_i]}\ \xi_{I^{'i}}.
\en

\end{cor}
In the trigonometric  and non-dynamical limit, the combinatorial structures of the formulas in this Corollary are the same as those in the 
geometric representation of $U_{q}(\slnh)$ on the equivariant $K$-theory of the quiver variety of type $A_{N-1}$  obtained by  Ginzburg and Vasserot\cite{GV, Vasserot}, and by Nakajima\cite{Na00}.  One can  directly check  that these actions of the  elliptic currents satisfy the defining relations of  the level-0 $U_{q,p}(\slnh)$ in the same way as in \cite{Vasserot,Na00}. We give a check of the most non-trivial relation \eqref{EFFE} in Appendix \ref{proofEFFE}.

\begin{prop}\lb{hwvec}
The finite dimensional representation given in Corollary \ref{actUqp} is an irreducible  
 highest weight representation with the highest weight vector $\xi_{11\cdots1}$. The elliptic analogue of the Drinfeld polynomials of this representation are given by
\bea
&& P_1(w)=\prod_{a=1}^n[u_a-v+1],\qquad P_l(w)=1\quad (l=2,\cdots,N-1) .
\ena
\noindent
{\it Proof.} The statement follows from a similar argument to  Theorem 4.11 in \cite{Konno09} and 
\be
&&E_j(1/w)\xi_{11\cdots1}=0\qquad (j=1,\cdots,N-1),\\
&&H^+_1(1/w)\xi_{11\cdots1}=\varrho \prod_{a=1}^n\frac{[u_a-v+1]}{[u_a-v]}\xi_{11\cdots1},\\
&&H^+_l(1/w)\xi_{11\cdots1}=\xi_{11\cdots1}\quad (l=2,\cdots,N-1) .
\en
\end{prop}

\section{Geometric Representation}\lb{GeomRep}
\subsection{Equivariant elliptic cohomology $\Ell_T(X)$}\lb{setting}
For $\la=(\la_1,\cdots,\la_N)\in \N^N, |\la|=n$, let $\F_\la$ denote the partial flag variety as before and consider the cotangent 
bundle $X=T^*\F_\la$. 
 Let us set  $T=A\times \C^*$, $A=(\C^*)^n$.  The torus $A$ has a natural action on $\F_\la$ and the extra $\C^*$ 
 acts on the fibers of $T^*\F_\la \to \F_\la$ by multiplication with weight  $\hbar$.  
Let 
$E=\C^*/p^{\Z}\ (|p|<1)$. We regard  the elliptic curve $E$ as a group scheme over $\C$.   
We follow \cite{AO,
FRV, 
Ga,Gro, GKV, RTV17
} for the definition of the $T$-equivariant 
elliptic cohomology $\Ell_T(X)$. 
 The basic facts on $\Ell_T(X)$ are summarized as follows.
\begin{itemize}
\item[(1)] 
The $T$-equivariant elliptic cohomology, $\Ell_T(X)$,  is a functor from finite $T$-spaces $X$ to superschemes,     
covariant in both $T$ and $X$, satisfying a set of axioms  (\cite{Gro}, 4.1 in \cite{Ga} and 2.1.2 in \cite{AO} ).  
In particular,  $\Ell_T(X)$ is  a scheme 
over $\Ell_T(\mathrm{ pt})\cong E^n\times E$. Moreover associated with  a construction of $X$ as a hyper-K\" ahler quotient
\be
&&T^*\PP\left(\bigoplus_{l=1}^{N-1}\Hom(\C^{\la^{(l)}},\C^{\la^{(l+1)}})\right)//G 
\en
with $G=\prod_{l=1}^{N-1}\mathrm{GL}({\la^{(l)}},\C)$, we have  a collection of   tautological vector bundles $\{\C^{\la^{(l)}}\}$ of 
$\mathrm{rk}=\la^{(l)}$ $(l=1,\cdots,N-1)$ over $X$ and a map
\bea
&&\Ell_T(X)\to\Ell_T({\rm pt})\times E^{(\la^{(1)})}\times \cdots \times E^{(\la^{(N-1)})},  \lb{fiberprod}
\ena
where $E^{(m)}=E^m/\gS_m$ denotes the symmetric product of $E$.
This map is expected to be an embedding near the origin of $\Ell_T(\pt)$ (2.2 in \cite{AO}).

\item[(2)] The Thom class map $ \Theta : \mathrm{K}_T(X)\to \Pic(\Ell_T(X))$ is a map of 
 a $T$-equivariant complex vector bundle $\xi$ to a line bundle $\mathbb{L}_T^\xi$ over $\Ell_T(X)$. The line bundle 
 $\mathbb{L}_T^\xi$ is called the Thom sheaf of $\xi$. See 2.3.2 in \cite{AO} and Definition 6.1 in \cite{Ga}.

\item[(3)] Let $f: X\to Y$ be a holomorphic map of $T$-spaces. Pull-back in the elliptic cohomology is the contravariant functoriality 
map $\Ell(f): \Ell_T(X) \to \Ell_T(Y)$ ( (1.7.4) in \cite{GKV} and 2.3.1 in \cite{AO}). 
If $f$ is proper, pushforward is a morphism 
\be
&&f_* : \Ell(f)_*\Theta(-N_f) \to \cO_{\Ell_T(Y)}
\en
of sheaves on $\Ell_T(Y)$, where $N_f=f^*TY-TX\in \mathrm{K}_T(X)$. 
See (2.3.2) in \cite{GKV} and (11) in \cite{AO}. 
\item[(4)] The dynamical parameter dependence is introduced by extending $\Ell_T(X)$
 to 
 \be
&&{\E}_T(X):=\Ell_T(X)\times {\cE}_{\Pic_T(X) }
\en
where 
\be
{\cE}_{\Pic_T(X)}=\Pic_T(X)\otimes_{\Z} E,
\en
as a scheme over $\B_{T, X}={\Ell}_T({\rm pt})\times {\cE}_{\Pic_T(X)}$. 
The variables in the two factors of $\B_{T, X}$, $z_1,\cdots,z_n$, $\hbar$ in ${\Ell}_T({\rm pt})\cong E^n\times E$ and  $\Pi^*_{j,j+1}\ (1\leq j\leq N-1)$ 
in $ {\cE}_{\Pic_T(X)}$, 
 are called the equivariant and the K\" ahler parameters, respectively. See 2.4.3 in \cite{AO}. 
In Sec.\ref{directCom}, we identify $\Pi^*_{j,j+1}$ with the dynamical parameters. 
\item[(5)] $\cE_{\Pic_T(X)}$ and $\cE^\vee_{\Pic_T(X)}:=\Hom(\Pic_T(X),E)$ are dual abelian varieties each other. 
Hence there exists a universal line bundle $\cU_{\mbox{\footnotesize Poincar\' e}}$ over $\cE^\vee_{\Pic_T(X)}\times \cE_{\Pic_T(X)}$. 
By using  a map 
\be
&&\tilde{c} : \Ell_T(X) \to \cE^\vee_{\Pic_T(X)}, 
\en
which is obtained from the Chern class of line bundles over $X$ (2.4.1 in \cite{AO}), 
one obtains  a line bundle $\cU$ on $\E_T(X)$ as
\be
&&\cU=(\tilde{c}\times 1)^*\ \cU_{\mbox{\footnotesize Poincar\' e}}.
\en
\end{itemize}

\subsection{Elliptic stable envelopes}

\subsubsection{Chamber structure}
Let  $\Hom_{\rm grp}(\C^*,A)$ be the space of one parameter subgroups $\rho$ in $A$ and  
$\Hom_{\rm grp}(\C^*,A)\otimes_\Z\R\subset \Lie A$ be its real form.  
 The latter space can be decomposed into finitely many  chambers $\gC$ 
defined as a connected component of the compliment of the union of hyperplanes given by $\rho$ such that $X^{\rho(\C^*)}\not=X^A$\cite{Na16}.  

The $A$-fixed points on $X$ are described by the partitions in $ \cI_\la$. 
Let $X^A$ be the $A$-fixed point locus in $X$ and  $X^A=\bigsqcup_{I\in \cI_\la} F_I$ a decomposition to connected components. 
Let $\rho\in \gC$. For every $\cS\subset X^A$ we define its attracting set 
\be
&&\mathrm{ Attr}(\cS)=\{(x,s), s\in \cS,  \lim_{t\to 0}\rho(t)x=s\ \}\ \subset\ X\times X^A, 
\en
and denote by $\mathrm{ Attr}^f(\cS)$ the full attracting set, which is the minimal closed subset of $X$ that contains the diagonal 
$\cS\times \cS$  and is closed under taking $\mathrm{ Attr}(\cdot)$.  We then define a partial ordering on $\{F_I\}$ by 
\be
&&F_J\leq F_I\quad \Leftrightarrow\quad \mathrm{ Attr}^f(F_I)\cap F_J\not=\emptyset.
\en 

\subsubsection{Definition}
For a pair $(\mu,\nu)$, $\mu\in \mathrm{char}(T)=\Hom(\Ell_T(\mathrm{pt}),E)$, $\nu\in \Pic_T(X)=\Hom(E,\cE_{\Pic_T(X)})$,  let $\varsigma$ denote the automorphism of $\cB_{T,X}$
\be
&&\varsigma(\mu\nu) : (z,\Pi^*)\mapsto (z,\nu(\mu(z)) \Pi^*).  
\en

Let $\iota :X^A\to X$ be the inclusion map, which is proper. 
For each chamber $\mathfrak{C}$ of $\mathrm{Lie}\, A$, one can consider  the polarization  $T^{1/2}X\in \mathrm{K}_T(X)$ of $X$
and its restriction $T^{1/2}X|_{X^A}$ to $X^A$.  Let us denote by $\mathrm{ind}:=T^{1/2}X|_{X^A,>0}$  the attracting part of $T^{1/2}X|_{X^A}$. We have 
 \be
 &&\mathrm{det}\, \ind
  \in \Pic_T(X^A) 
 \en
and a translation
\be
&&\varsigma(-\hbar\, \mathrm{det}\, \ind
) : \cB_{T,X^A} \to \cB_{T,X^A}. 
\en
For  the line bundle $\cU_{\Ell_T(X^A)}$ on 
 $\Ell_T(X^A)\times \cE_{\Pic_T(X)}$  we set
 \be
 &&\cU'=(1\times \iota^*)^*\varsigma(-\hbar\, \mathrm{det}\, \ind
 )^* \, \cU_{\Ell_T(X^A)},
 \en
 where $\iota^*$ is the pull-back of line bundles from $X$ to $X^A$.

The elliptic stable envelop ${\rm Stab_{\mathfrak{C}}}$ is defined to be a map of $\cO_{\B_{T,X}}$-modules
\be
&&
\Theta(T^{1/2}X^A)\otimes  \cU' \to 
 \Theta(T^{1/2}X)\otimes \cU \otimes \cdots  ,
\en 
where $\Theta(T^{1/2}X)$ denotes the Thom sheaf of a polarization, 
and $\cdots$ stands for 
a certain line bundle pulled back from
\be
&&\B'=\B_{T,X}/\Ell_A(\mathrm{pt}).
\en 
${\rm Stab_{\mathfrak{C}}}$ is subjected to the following two conditions (3.3.4 in \cite{AO}).
\begin{itemize}
\item[(i)] (triangularity) Let $s_K$ be an elliptic cohomology class supported on $F_K$ locally over $\cB_{T,X}$. 
Then $\Stab_\gC(s_K)$ is supported on $\mathrm{Attr}^f(F_K)$. In particular if $F_K < F_I$ we have 
\be
&&\Stab_\gC(s_K)|_{F_I}=0.
\en
\item[(ii)] (normalization) Near the diagonal in $X\times F_K$, we have 
\be
&&\Stab_\gC =(-1)^{\mathrm{rk}\ \ind}\ j_*\pi^*,
\en
where
\be
&&F_K\stackrel{\pi}{\longleftarrow }\Attr(F_K)\stackrel{j}{\longrightarrow}X
\en 
are the natural projection and inclusion maps. 
\end{itemize}

\subsection{Direct comparison with the elliptic weight functions}\lb{directCom}
In this and the following subsections, we consider the elliptic weight functions  $\cW_I(t,z,\Pi^*)$ obtained from 
$\cW_I(t,z,\Pi)$ in \eqref{def:cW} by replacing $\Pi_{j,l}=q^{2(P+h)_{j,l}}$ by $\Pi^*_{j,l}=q^{2P_{j,l}}$, which also satisfy 
all the properties in Sec. \ref{EWF} under the same replacement. 

The symmetry structure in the target of \eqref{fiberprod} coincides with the one of the elliptic weight function 
$\cW_I(t,z,\Pi^*)$  with respect to   
the variables $t^{(l)}_a\ (l=1,\cdots,N-1, a=1,\cdots,\la^{(l)})$. This suggests that  $\{t^{(l)}_a\}$ can be identified with 
 the Chern roots of the tautological vector bundles over $X$\cite{BCFK}.  
 This structure as well as  the quasi periodicity  in Proposition \ref{quasiperiod} allow us to identify the elliptic weight functions $ \cW_I(t,z,\Pi^*)$ with meromorphic sections of line bundles over $\E_T(X)$ near the origin of  $\cB_{T,X}$. 

More precicely 
one can compare the elliptic weight functions $\cW_{\sigma_0(I)}(\tit,\sigma_0(z^{-1}),\Pi^{*-1})$, where $\tit$ denotes a set of $t$-variables associated with the partition $\sigma_0(I)$,  with the abelianization formula of the elliptic stable envelopes\cite{Shenfeld,AO}. 

Let $G$ be a reductive group acting on a vector space $M$. It induces the hamiltonian action on $T^*M$. 
Let $\mu_G$ be the corresponding moment map. Let  $S\subset G$ be the maximal torus and let $\pi_S:(\Lie G)^*\to (\Lie S)^*$
 be a projection, and set $\mu_S=\pi_S\circ \mu_G$. For the hyper-K\" aler quotient 
\be
&&X=\mu_G^{-1}(0)//G, 
\en
the associated abelian quotient
\be
&&X_S=\mu_S^{-1}(0)//S
\en
 is a hypertoric variety called the abelianization of $X$.

 Let $\la=(\la_1,\cdots,\la_N)\in \N^N$, $I\in \cI_\la$ and $\la^{(l)}=\la_1+\cdots+\la_l, I^{(l)}=\{i^{(l)}_1<\cdots<i^{(l)}_{\la^{(l)}}\}\ (l=1,\cdots,N)$ as in Sec.\ref{combNot}.

\subsubsection{The case $X=T^*\PP(\C^{(\la^{(l+1)})})$}\lb{TsP}
We follow a construction of $\Stab_\gC$ for $X$ given in 3.4  from \cite{AO}.  
For $a\in [1,\la^{(l)}]$, let us fix  a variable $t^{(l)}_a$ and regard it as 
the coordinate on $GL(1)/p^{\Z}=E=E^\vee$. The variable $t^{(l)}_a$  gives the Chern root of the line bundle $\cO(1)$ over $X$. 
 We denote by 
$\Pi^*_{l,l+1}=q^{2P_{l,l+1}}$ the K\" ahler parameter dual to $t^{(l)}_a$.

Let us  take a basis $\{e_{i^{(l+1)}_b}\ (b=1,\cdots,\la^{(l+1)})\}$ of $\C^{\la^{(l+1)}}$, on which  
the torus $A^{(l+1)}$ acts as $\mathrm{diag}({t}^{(l+1)}_1,\cdots,{t}^{(l+1)}_{\la^{(l+1)}})$. Hence   
 $F_{i_b^{(l+1)}}=\C e_{i^{(l+1)}_b}$ $(b=1,\cdots,\la^{(l+1)})$ give the $A^{(l+1)}$-fixed points.
We chose  the chamber $\gC^{(l+1)}$  such that 
\be
&&F_{i^{(l+1)}_1}>F_{i^{(l+1)}_2}>\cdots>F_{i^{(l+1)}_{\la^{(l+1)}}}. 
\en

Let  $\phi^{(l)}: [1,\la^{(l)}] \to [1,\la^{(l+1)}]$ be a map defined by $i^{(l+1)}_{\phi^{(l)}(a)}=i^{(l)}_a$. 
For $a\in [1,\la^{(l)}]$, the elliptic stable envelopes $\Stab_{\gC^{(l+1)}}(F_{\phi^{(l)}(a)})$ for 
$X=T^*\PP(\C^{\la^{(l+1)}})$ is then given by 
\bea
&&f^{(l)}(t^{(l)}_a,t^{(l+1)},\Pi^*_{l,l+1}q^{2(\la^{(l+1)}-\phi^{(l)}(a))})\nn\\
&&\qquad=
\left.
\frac{[-v^{(l+1)}_{\phi^{(l)}(a)}+v^{(l)}_a+P_{l,l+1}+\la^{(l+1)}-\phi^{(l)}(a)
]}{[P_{l,l+1}+\la^{(l+1)}-\phi^{(l)}(a)
]}\right.
\nn\\
&&\qquad\qquad\times 
\prod_{b=1\atop b<\phi^{(l)}(a)}^{\la^{(l+1)}}[v^{(l+1)}_b-v^{(l)}_a]
\prod_{b=1\atop b>\phi^{(l)}(a)}^{\la^{(l+1)}}[v^{(l+1)}_b-v^{(l)}_a+1], \lb{def:fl}
\ena
where we set  $t^{(l)}_a=q^{2v^{(l)}_a}$, $t^{(l+1)}_a=q^{2v^{(l+1)}_a}$, 
 $t^{(l+1)}=(t^{(l+1)}_1,\cdots,t^{(l+1)}_{\la^{(l+1)}})$, and made an identification  $\hbar=q^{-2}$.  
 
Let us compare this with  \eqref{def:u}. 
Let $\tilde{I}=\sigma_0(I)$  with $\tilde{I}^{(l)}=\{\ti^{(l)}_1<\cdots <\ti^{(l)}_{\la^{(l)}}\}$ $(l=1,\cdots,N)$ and consider the variables $\tit^{(l)}_{a}\equiv t(\ti^{(l)}_a)$. 
Then one finds 
\bea
&&\ti^{(l)}_a=\sigma_0\left(i^{(l)}_{\sigma_0^{(l)}(a)}\right), \lb{rel:ipi}
\ena 
where $\sigma_0^{(l)}\in \gS_{\la^{(l)}}$
denotes the longest element.
Furthermore for a map $\tilde{\phi}^{(l)}: [1,\la^{(l)}] \to [1,\la^{(l+1)}]$ such that  $\ti^{(l+1)}_{\tilde{\phi}^{(l)}(a)}=\ti^{(l)}_a$, 
we have 
\bea
&&\tilde{\phi}^{(l)}(\sigma^{(l)}_0(a))=\sigma_0^{(l+1)}(\phi^{(l)}(a)). \lb{rel:mupmu}
\ena
We then find  that \eqref{def:fl}  coincides with 
\be
&&u^{(l)}_{\sigma_0(I)}(\tit^{(l)}_{\ta}, \tit^{(l+1)},\Pi^{*-1}_{l,l+1}q^{-2(\tilde{\phi}^{(l)}(\ta)-1)})
\en 
from \eqref{def:u} under the identification  
 $\ta=\sigma^{(l)}_0(a)$,  $\tb=\sigma^{(l+1)}_0(b)$, $t^{(l)}_a=t(i^{(l)}_a)$, $t^{(l+1)}_b=t(i^{(l+1)}_b)$, 
$t^{(l)}_a=\tit^{(l)}_{\ta}$, $t^{(l+1)}_b=\tit^{(l+1)}_{\tb}$ and in particular $t^{(N)}_a=z^{-1}_{\sigma_0(a)}$. This is due to the  identity  
\bea
&&\tilde{\phi}^{(l)}(\ta)-1=\la^{(l+1)}-\phi^{(l)}(a),
\ena
and the equality such as 
\bea 
&& \prod_{\tb=1\atop \ti^{(l+1)}_{\tb}<\ti^{(l)}_{\ta}}^{\la^{(l+1)}}[{\tv}^{(l+1)}_{\tb}-{\tv}^{(l)}_{\ta}+1]=
\prod_{b=1\atop b>\phi^{(l)}(a)}^{\la^{(l+1)}}[v^{(l+1)}_b-v^{(l)}_a+1], 
\ena
which follows from  
\be
&& \ti^{(l+1)}_{\tb}<\ti^{(l)}_{\ta}\quad  \Leftrightarrow\quad  \tb<\tilde{\phi}^{(l)}(\ta)\quad \Leftrightarrow\quad b>\phi^{(l)}(a). 
\en

\subsubsection{The case  $X=T^*\mathrm{Gr}(\la^{(l)},\la^{(l+1)})$}
We follow 4.4  in \cite{AO}. 
For each $l$ $(l=1,\cdots,N-1)$, let  $S^{(l)}=\mathrm{diag}({t}^{(l)}_1,\cdots,{t}^{(l)}_{\la^{(l)}})\subset GL(\la^{(l)})$ be the maximal torus. 
The abelianization  $X_{S^{(l)}}$ of $X$ is given by 
\be
&&X_{S^{(l)}}=(T^*\PP(\C^{\la^{(l+1)}}))^{\la^{(l)}}. 
\en
Correspondingly the product 
\be
&&\Stab^{S^{(l)}}_{\gC^{(l+1)}}(F_{I^{(l)}})=\prod_{a=1}^{\la^{(l)}}f^{(l)}(t^{(l)}_a,t^{(l+1)},\Pi^*_{l,l+1}
q^{2C_{l,l+1}(i^{(l)}_a)})
\en
gives the elliptic stable envelopes for $X_{S^{(l)}}$. Here $F_{I^{(l)}}=\mathrm{Span}_\C\{\ e_{i^{(l+1)}_b} \ |\ i^{(l+1)}_b=i^{(l)}_a\ ( a=1,\cdots,\la^{(l)})\}=
\mathrm{Span}_\C\{\ e_{i^{(l)}_a} \ ( a=1,\cdots,\la^{(l)})\}$  are the 
$A^{(l+1)}$-fixed points, and ${t}^{(l)}_a$'s are identified with the Chern roots of the tautological bundle over
$X$. The dynamical shift $C_{l,l+1}(i^{(l)}_a)$ is given by 
\bea
&&C_{l,l+1}(i^{(l)}_a)=2(\la^{(l)}-a)-\la^{(l+1)}+\phi^{(l)}(a). \lb{Cllp1}
\ena
Note that this is identical to the one  in Proposition \ref{combC}  for $i^{(l)}_a=s, \mu_s=l$.

Then  the  abelianization formula\cite{Shenfeld,AO} gives 
the stable envelopes $\Stab_{\gC^{(l+1)}}(F_{I^{(l)}})$ for  $X$ in terms of  $\Stab^{S^{(l)}}_{\gC^{(l+1)}}(F_{I^{(l)}})$ 
 as follows. 
\be
&&\Stab_{\gC^{(l+1)}}(F_{I^{(l)}})={\rm Sym}_{t^{(l)}} \frac{\Stab^{S^{(l)}}_{\gC^{(l+1)}}(F_{I^{(l)}})
}{\prod_{1\leq a<b\leq \la^{(l)}}[v^{(l)}_a-v^{(l)}_b][v^{(l)}_b-v^{(l)}_a-1]}.
\en

In particular, the case $N=2$, where  $\la^{(2)}=n$ and $t^{(2)}_b=z^{-1}_{b}$, we have  for $a\in [1,\la^{(1)}]$ 
\be
&&f^{(1)}(t^{(1)}_a,z^{-1},\Pi^*_{1,2}q^{2C_{1,2}(i^{(1)}_a)})\\
&&=\left.
\frac{[u_s+v^{(1)}_a+P_{1,2}+C_{1,2}(i^{(1)}_a)
]}{[P_{1,2}+C_{1,2}(i^{(1)}_a)
]}\right.
\prod_{1\leq b<\phi^{(1)}(a)}[u_{b}+v^{(1)}_a]
\prod_{\phi^{(1)}(a)<b\leq n}[u_{b}+v^{(1)}_a-1].
\en 
The resultant   $\Stab_{\gC^{(2)}}(F_{I^{(1)}})$ is the elliptic stable envelope for $T^*\Gr(\la^{(1)},n)$ given in  (60) from \cite{AO}. 
Combining this with the identification between $f^{(l)}$ and $u^{(l)}_{\sigma_0(I)}$ for $l=1$ in the last subsection, we  find 
\be
&&\Stab_{\gC^{(2)}}(F_{I^{(1)}})=\cW_{\sigma_0(I)}(\tit,\sigma_0(z^{-1}),\Pi^{*-1}),
\en 
where 
$\tit=(\tit^{(l)}_{a})\ (l=1,2,\ a=1,\cdots,\la^{(1)})$.

\subsubsection{The case  $X=T^*\F_\la$}
Let $N>2$. The partial flag variety $\F_\la=\F(\la^{(1)},\cdots,\la^{(N-1)},n)$ is given by a hyper-K\" ahler quotient 
\be
&&\F_\la=\PP\left(\bigoplus_{l=1}^{N-1}\Hom(\C^{\la^{(l)}},\C^{\la^{(l+1)}})\right)//G
\en
by the action of $G=\prod_{l=1}^{N-1}GL(\la^{(l)},\C)$, and has a description as a tower of Grassmannian bundles\cite{BCFK}
\be
 \Gr(\la^{(l)},\la^{(l+1)})\quad &\to&\quad \F(\la^{(l)},\la^{(l+1)}, \cdots,\la^{(N-1)},n) =\Gr(\la^{(l)},\cS_{l+1})\\
&&\hspace{2cm}\downarrow\\
&&\quad \F(\la^{(l+1)}, \cdots,\la^{(N-1)},n),
\en
where $\cS_1\subset \cS_2\subset \cdots \subset \cS_{N-1}\subset \C^n\otimes \cO_{\F}$ denotes the universal bundle. 
Correspondingly, the abelianization $(\F_{\la})_S$ of $\F_\la$ by $S=\prod_{l=1}^{N-1}S^{(l)}\subset G$ is a tower of product-of-projective-space bundles:
\be
(\PP^{\la^{(l+1)}-1})^{\la^{(l)}}\quad &\to&\quad \F_l\quad  =\quad \PP(V_{l+1})\times_{\F_{l+1}}\cdots \times_{\F_{l+1}} \PP(V_{l+1})\\
&&\quad\ \downarrow\\
&&\quad \F_{l+1},
\en  
where $\F_l$ is the abelianization of $\F(\la^{(l)},\la^{(l+1)}, \cdots,\la^{(N-1)},n)$, and 
\be
&&V_{l+1}=\bigoplus_{j=1}^{\la^{(l+1)}}\cO(0,\cdots,0,\stackrel{j}{-1},0,\cdots,0). 
\en
is the vector bundle on $\F_{l+1}$ corresponding to $\cS_{l+1}$. 
The abelianization $X_S$ of $X=T^*\F_\la$ containes $T^*(\F_\la)_S$ as a dense open subset, and  is
 obtained by gluing together the flopped versions of $T^*(\F_\la)_S$ (5.1.2 in \cite{Shenfeld}). 

For $I=I_{\mu_1,\cdots,\mu_n}\in \cI_\la$, let $F_I$ denote the $A=\mathrm{diag}(z_1,\cdots,z_n)$-fixed point 
$0\subset F_{I^{(1)}}\subset \cdots \subset F_{I^{(N-1)}}\subset \C^{n}$ and $\{t^{(l)}_a\}$ 
be the Chern roots of the tautological bundle over $X$.   We choose the chamber 
$\gC$  consistently to the choice of  $\gC^{(l+1)}$ for $T^*\Gr(\la^{(l)},\la^{(l+1)})$ $(l=1,\cdots,N-1)$. 
Namely  
\be
&&z_b/z_a>0 \quad \Leftrightarrow\quad a<b.
\en

Then the  elliptic stable envelop for $X_S$ is given by
\bea
\Stab^S_{\gC}(F_I)&=&\prod_{l=1}^{N-1}
\prod_{a=1}^{\la^{(l)}}f^{(l)}\left(t^{(l)}_a,t^{(l+1)},\Pi^*_{\mu_{i^{(l)}_a},l+1}
q^{2C_{\mu_{\mbox{\tiny$i^{(l)}_a$}},l+1}(i^{(l)}_a)}\right).\lb{StabXS}
\ena
Here we  defined  for $s\in I^{(l)}$, $\mu_s\leq l$,  
\bea
&&\Pi^*_{\mu_{s},l+1}:=\prod_{k=\mu_{s}}^l \Pi^*_{k,k+1}, \qquad\qquad
q^{2C_{\mu_s,l+1}(s)}:=\prod_{k=\mu_s}^l q^{2C_{k,k+1}(s)} \lb{factorization}
\ena
with $C_{k,k+1}(s)$ given in \eqref{Cllp1} for $s=i^{(l)}_a$.  
The factor $\ds{\prod_{a=1}^{\la^{(l)}}f^{(l)}\left(t^{(l)}_a,t^{(l+1)},\Pi^*_{\mu_{i^{(l)}_a},l+1}
q^{2C_{\mu_{i^{(l)}_a},l+1}(i^{(l)}_a)}\right)}$ is 
 a contribution from $T^*\Gr(\la^{(l)},\la^{(l+1)})$ with the  K\" ahler parameters $\Pi^*_{\mu_{i^{(l)}_a},l+1}$  
 dual to $t^{(l)}_a$. Note that  for  $N>2$, $\mu_{i^{(l)}_a}$  in general takes a value in the range $[1,l]$ corresponding  to the embedding structure $F_{I^{(1)}}\subset \cdots\subset F_{I^{(\mu_{i^{(l)}_a})}}\subset \cdots \subset F_{I^{(l)}}\subset F_{I^{(l+1)}}$.   
The formulas in  \eqref{factorization} are then natural in the sense  that the  K\" ahler parameters as well as their dynamical shift 
$q^{2C_{\mu_{i^{(l)}_a},l+1}(i^{(l)}_a)}$ are given as built-up contributions 
from a sequence of the  Grassmannians filling a gap between $F_{I^{(\mu_{i^{(l)}_a})}}$ and $F_{I^{(l+1)}}$. 
 Note that  the resultant 
  $C_{\mu_s,l+1}(s)$ coincides with the one in Proposition \ref{combC}.

Then again the abelianization formula\cite{Shenfeld, AO} yields the following expression of 
$\Stab_{\gC}(F_I)$ for $X=T^*\F_\la$. 
\bea
\Stab_{\gC}(F_I)&=&{\rm Sym}_{t^{(1)}}\cdots {\rm Sym}_{t^{(N-1)}}
\frac{\Stab^S_{\gC}(F_I)
}{\prod_{l=1}^{N-1}\prod_{1\leq a<b\leq \la^{(l)}}[v^{(l)}_a-v^{(l)}_b][v^{(l)}_b-v^{(l)}_a-1]}\nn\\
&=&{\rm Sym}_{t^{(1)}}\cdots {\rm Sym}_{t^{(N-1)}}\prod_{l=1}^{N-1}
\frac{
\prod_{a=1}^{\la^{(l)}}f^{(l)}(t^{(l)}_a,t^{(l+1)},\Pi^*_{\mu_{i^{(l)}_a},l+1}q^{2C_{\mu_{i^{(l)}_a},l+1}({i^{(l)}_a})})
}
{\prod_{1\leq a<b\leq \la^{(l)}}[v^{(l)}_a-v^{(l)}_b][v^{(l)}_b-v^{(l)}_a-1]}.\nn\\
&&\lb{StabPF}
\ena  
Under the same identification as in  the previous subsections, we thus obtain 
\bea 
&&\Stab_{\gC}(F_I)= \cW_{\sigma_0(I)}(\tit,\sigma_0(z^{-1}),\Pi^{*-1}). \lb{StabtW}
\ena

\subsubsection{Restriction to the fixed points}
For $J\in \cI_\la$, let us consider the restriction to the fixed point $F_J$ i.e. the specialization $t=z^{-1}_{J}$  given by 
\be
&&{t}^{(l)}_a=z^{-1}_{j^{(l)}_a} \qquad (l=1,\cdots,N-1, a=1,\cdots,\la^{(l)} ). 
\en
From \eqref{StabPF} we obtain
\be
&&\Stab_\gC(F_I)|_{F_J}={\rm Sym}_{t^{(1)}}\cdots {\rm Sym}_{t^{(N-1)}}\prod_{l=1}^{N-1}
\frac{
\prod_{a=1}^{\la^{(l)}}f^{(l)}(t^{(l)}_a,t^{(l+1)},\Pi^*_{\mu_{i^{(l)}_a},l+1}q^{2C_{\mu_{i^{(l)}_a},l+1}({i^{(l)}_a})})|_{t=z^{-1}_{J}}
}{\prod_{1\leq a<b\leq \la^{(l)}}[u_{j^{(l)}_a}-u_{j^{(l)}_b}][u_{j^{(l)}_b}-u_{j^{(l)}_a}+1]},
\en
where 
\be
&&f^{(l)}(t^{(l)}_a,t^{(l+1)},\Pi^*_{\mu_{s},l+1}q^{2C_{\mu_s,l+1}(s)})|_{t=z^{-1}_{J}}\\
&&\qquad=
\left.
\frac{[u_{j^{(l+1)}_b}-u_{j^{(l)}_a}+P_{\mu_s,l+1}+C_{\mu_s,l+1}(s)]
}{[P_{\mu_s,l+1}+C_{\mu_s,l+1}(s)
]}\right|_{i^{(l+1)}_b=i^{(l)}_a=s}
\\
&&\qquad\times 
\prod_{b=1\atop i^{(l+1)}_b<s}^{\la^{(l+1)}}[u_{j^{(l+1)}_b}-u_{j^{(l)}_a}]
\prod_{b=1\atop i^{(l+1)}_b>s}^{\la^{(l+1)}}[u_{j^{(l+1)}_b}-u_{j^{(l)}_a}-1].
\en
This is an elliptic and dynamical analogue of the formula in Theorem 5.2.1 in \cite{Shenfeld}. 

By using  the identification \eqref{StabtW} and the equivalence of the specializations
\be
&&t=z^{-1}_J \quad \Leftrightarrow\quad \tit=\sigma_0(z^{-1})_{\sigma_0(J)}, \quad \mbox{i.e.} \quad 
t^{(l)}_a=z^{-1}_{j^{(l)}_a} \quad \Leftrightarrow\quad \tit^{(l)}_{\ta}=z^{-1}_{\sigma_0(\tj^{(l)}_{\ta})},
\en
where $\tilde{J}=\sigma_0(J)$, we obtain the following  identification.
\bea
 {\rm Stab}_\gC(F_I)\vert_{F_J}&=&\cW_{\sigma_0(I)}(\sigma_0(z^{-1})_{\sigma_0(J)},\sigma_{0}(z^{-1}),\Pi^{*-1}) \nn\\
&=&\cW_{\sigma_0(I)}(z^{-1}_{J},\sigma_{0}(z^{-1}), \Pi^{*-1}), \lb{StabWFP}
\ena
where in the second equality we used the identity
\be
&&\cW_{I}(\sigma(z)_{\sigma(J)},z,\Pi^*)=\cW_{I}(z_{J},z,\Pi^*)\qquad \forall I,J \in \cI_\la, \ \forall \sigma\in \mathfrak{S}_n. 
\en

\subsection{Geometric representation}\lb{subGeomRep}
Let $X=T^*\F_\la$ and fix a chamber $\mathfrak{C}$ as above. 
By definition,  the stable classes ${\rm Stab}_\gC(F_K)$ $(K\in \cI_\la)$  are triangular with respect to the fixed point classes $\{[I]\}_{I\in \cI_\la}$ 
 in $\E_T(X)$.   See 3.3.4 in \cite{AO}. Namely we have the following expansion formula
\bea
&&{\rm Stab}_\gC(F_K)=\sum_{I\in \cI_\la}\frac{{\rm Stab}_\gC(F_K)\vert_{F_I}}{R(z^{-1}_I)}\ [I].\lb{StabI}
\ena
Here we  chose a normalization by $R(z_I)$ given 
 in Proposition \ref{orthogonalProp} for later convenience. We regard this as the definition of the fixed point classes. 

From \eqref{StabWFP}, we have 
\bea
&&{\rm Stab}_\gC(F_K)\vert_{F_I}
=\cW_{\sigma_0(K)}(z^{-1}_{I},\sigma_{0}(z^{-1}),\Pi^{*-1}). \lb{StabW}
\ena 
Note also that by the replacement $z\mapsto z^{-1}$ and $\Pi \mapsto \Pi^{*-1}$ one can rewrite Proposition \ref{orthogonalProp} 
 as 
\be
&&\sum_{I\in \cI_\la}\frac{\cW_J(z^{-1}_I,z^{-1},\Pi^* q^{2\sum_{j=1}^{n}<\bep_{\mu_j},h>})
\cW_{\sigma_0(K)}(z^{-1}_I,\sigma_0(z^{-1}),\Pi^{*-1})}{Q(z^{-1}_I)R(z^{-1}_I)}=\delta_{J,K}.
\en
Then using this and \eqref{StabW} one can invert \eqref{StabI} and obtain  
\bea
[I]=\sum_{J\in \cI_\la} \tW_{J}(z^{-1}_I,z^{-1},\Pi^*  q^{2\sum_{j=1}^n<\bep_{\mu_j},h>})\ {\rm Stab}_\gC(F_J). \lb{IStab}
\ena
Comparing this with 
Theorem \ref{XtW}, we find that 
\eqref{IStab} is identical to the relation \eqref{xiv} under the 
correspondence 
\be
 \mbox{the Gelfand-Tsetlin base}\ \xi_I &\Leftrightarrow& \mbox{the fixed point class }\ [I],\\
 \mbox{the standard base }\ v_J &\Leftrightarrow& \mbox{the stable class}\   {\rm Stab}_\gC(F_J). 
\en

On the basis of  this correspondence  as well as  Theorem \ref{actHC} and Corollary \ref{actUqp}, we obtain the following statement on the level-0 action of $E_{q,p}(\glnh)$ and $U_{q,p}(\slnh)$ on ${\E}_T(X)$.  

\begin{thm}\lb{geomactHC}
Under the same notation as Theorem \ref{actHC}, 
let us define the action of the half-currents $\cK_j^\pm(1/w), E^\pm_{j+1,j}(1/w,P), F^\pm_{j,j+1}(1/w,P)$ on the fixed point classes 
by
\bea
&&\cK_j^\pm(1/w)[I]=\left.\prod_{k=1}^{j-1}\prod_{a\in I_k}\frac{[u_a-v]}{[u_a-v+1]}\right|_{\pm}\left.\prod_{l=j+1}^N\prod_{b\in {I_l}}\frac{[u_b-v-1]}{[u_b-v]}\right|_{\pm}\ [I],\lb{GactK}\\
&&\hspace{-0.2cm}E^\pm_{j+1,j}(1/w,P)[I]=\sum_{i\in I_{j+1}}\left.\frac{[P_{j,j+1}-u_i+v][1]}{[P_{j,j+1}][u_i-v]}\right|_{\pm}\prod_{k\in I_{j+1}\atop \not=i}\frac{[u_i-u_k+1]}{[u_i-u_k]}\ [{I^{i'}}],\lb{GactE}\\
&&\hspace{-0.7cm}F^\pm_{j,j+1}(1/w,P)[I]=\sum_{i\in I_{j}}\left.\frac{[P_{j,j+1}+\la_j-\la_{j+1}+u_i-v-1][1]}{[P_{j,j+1}+\la_j-\la_{j+1}-1][u_i-v]}\right|_{\pm}\prod_{k\in I_{j}\atop \not=i}\frac{[u_k-u_i+1]}{[u_k-u_i]}\ [{I^{'i}}]. \lb{GactF}
\ena
Then this gives an irreducible finite-dimensional representation of $E_{q,p}(\glnh)$  on ${\E}_T(X)$. 
\end{thm}

\begin{cor}\lb{geomactUqp}
The level-0 action of  $U_{q,p}(\slnh)$ on ${\E}_T(X)$ is given by
\be
&&H_j^\pm(q^{j-N+1}/w)[I]=\varrho
\left.\prod_{a\in I_j}\frac{[u_a-v+1]}{[u_a-v]}\right|_{\pm}\left.\prod_{b\in {I_{j+1}}}\frac{[u_b-v-1]}{[u_b-v]}\right|_{\pm}\ e^{-Q_{\al_j}}[I],\\
&&E_{j}(q^{j-N+1}/w)[I]=\frac{\mu^*[1]}{[0]'}
\sum_{i\in I_{j+1}}\delta(z_i/w)\prod_{k\in I_{j+1}\atop \not=i}\frac{[u_i-u_k+1]}{[u_i-u_k]}\ e^{-Q_{\al_j}}[{I^{i'}}],\\
&&F_{j}(q^{j-N+1}/w)[I]=\frac{\mu[1]}{[0]'}
\sum_{i\in I_{j}}\delta(z_i/w)\prod_{k\in I_{j}\atop \not=i}\frac{[u_k-u_i+1]}{[u_k-u_i]}\ [{I^{'i}}].
\en
\end{cor}

\vspace{2mm}
\noindent
{\it Remark.} Similar correspondences between the Gelfand-Tsetlin basis and the fixed point classes were studied in \cite{Nagao, Kodera, FFFR, FFNR,Tsym}. In \cite{Nagao}  for the level-(0,1) representation of the quantum toroidal algebra of type $A$, the Gelfand-Tsetlin basis on the $q$-Fock space\cite{Uglov} was identified with the fixed point basis 
of the equivariant $\mathrm{K}$-theory of corresponding cyclic quiver variety\cite{VV99}. Affine Yangian analogue of this result was obtained in \cite{Kodera}.  In \cite{FFFR, FFNR,Tsym},  certain geometric actions of the universal enveloping  algebra $U(\gln)$  on the Laumon spaces, of the affine Yangian of type $A^{(1)}_{N-1}$  on the affine Laumon spaces and of the quantum toroidal algebra $\ddot{U}_{q}(\slnh)$ on the $\mathrm{K}$-theory of the affine Laumon spaces were constructed, respectively.

\section*{Acknowledgements}
The author would like to thank Andrei Okounkov for stimulating discussions. 
He is also grateful to 
Vassily Gorbounov, Michio Jimbo, Yoshiyuki Kimura, Christian Korff, 
Andrei Negut,  Masatoshi Noumi, Vitaly Tarasov, Yaping Yang and Paul Zinn-Justin  
for  discussions, and to the organizers of the workshops 
{\it ``Elliptic Hypergeometric Functions in Combinatorics,
Integrable Systems and Physics''} and {\it ``Geometric $R$-Matrices: from geometry to probability''} 
for kind invitation. 
His research is supported by the Grant-in -Aid for Scientific Research (C) 17K05195, JSPS.

\begin{appendix}

\section{The $H$-Hopf algebroid structure of $E_{q,p}(\glnh)$ and $U_{q,p}(\glnh)$}\lb{HHopfalgebroid}
We summarize  an $H$-Hopf algebroid structure based on the opposite coproduct to the one used in the previous papers\cite{Konno09, Konno16}. This opposite $H$-Hopf algebroid structure is used in  Sec.\ref{finRep} to construct the finite dimensional representations of $E_{q,p}(\glnh)$ and $U_{q,p}(\glnh)$. 

Let $\cA$ denote  $E_{q,p}(\glnh)$  or $U_{q,p}(\glnh)$.
The $\cA$ is  bi-graded over $H^*$  by
\be
&&\cA=\bigoplus_{\al,\beta\in H^*}\cA_{\al,\beta}\\
&&\cA_{\al,\beta}=\left\{x\in \cA \left|\ q^{P+h}x q^{-(P+h)}=q^{<\al,P+h>}x,\quad q^{P}x q^{-P}=q^{<\beta,P>}x\ \forall P+h, P\in H\right.\right\}, 
\en
and possesses two moment maps $\mu_l, \mu_r : \FF \to \cA_{0,0}$ defined by 
\be
&&\mu_l(\hf)=f(P+h,p)\in \FF[[p]],\qquad \mu_r(\hf)=f(P,p^*)\in \FF[[p]].
\en
The $\mu_l$ and $\mu_r$ satisfy
\be
\mu_l(\hf)a=a \mu_l(T_\al \hf), \quad \mu_r(\hf)a=a \mu_r(T_\beta \hf), \qquad 
a\in \cA_{\al,\beta},\ \hf\in \FF,
\en
where $T_\al=e^\al\in \C[\cR_Q]$ denotes the automorphism of $\FF$ 
\be
{T}_{\al}\widehat{f}=e^{\al}f(P,P+h)e^{-\al}={f}(P+<\al,P>,P+h+<\al,P>).
\en
Hence  both  $E_{q,p}(\glnh)$  and $U_{q,p}(\glnh)$ are the $H$-algebras\cite{EV,KR,Konno16}. 

Let $\cA$ and $\cB$ be two $H$-algebras. The tensor product $\cA {\widetilde{\otimes}}\cB$ is the $H^*$-bigraded vector space with 
\bea
 (\cA {\widetilde{\otimes}}\cB)_{\al\beta}=\bigoplus_{\gamma\in H^*} (\cB_{\gamma\beta}\otimes_{\FF}\cA_{\al\gamma}),\lb{def:tot1}
\ena
where $\otimes_{\FF}$ denotes the usual tensor product 
modulo the following 
relation.
\bea
\mu_l^\cB(\hf) b\otimes a=b\otimes \mu_r^\cA(\hf) a, \qquad a\in \cA, \ 
b\in \cB, \ \hf\in \FF.\lb{AtotB}
\ena
The tensor product $\cA {\widetilde{\otimes}}\cB$ is again an $H$-algebra with the multiplication $(b\otimes a)(d\otimes c)=bd\otimes ac$ 
$(a,c\in \cA, b,d\in \cB)$ and the moment maps 
\bea
\mu_l^{\cA {\widetilde{\otimes}}\cB} =1\otimes \mu_l^\cA,\qquad \mu_r^{\cA {\widetilde{\otimes}}\cB} =\mu_r^\cB\otimes 1.\lb{def:tot2}
\ena

We also consider the $H$-algebra of the shift operators\cite{EV} 
\be
&&\cD=\{\ \sum_\al \widehat{f}_\al{T}_{\al}\ |\ \widehat{f}_\al\in {M}_{H^*}, 
\al\in \cR_Q\ \},\\
&&\cD_{\ha,\ha}=\{\ \widehat{f}{T}_{-\ha}\ \},\quad \cD_{\ha,{\beta}}=0\ 
(\al\not=\beta),\\
&&\mu_l^{\cD}(\widehat{f})=
\mu_r^{\cD}(\widehat{f})=\widehat{f}{T}_0 \qquad \widehat{f}\in {M}_{H^*}.
\en
Then we have the $H$-algebra isomorphism 
\bea
 \cA\cong \cA\tot\cD\cong \cD\tot \cA. \lb{Diso}
\ena 

The two $H$-algebras  $E_{q,p}(\glnh)$ and  $U_{q,p}(\glnh)$ are equipped with  the common $H$-Hopf algebroid structure\cite{Konno16}  
 defined  by the two $H$-algebra homomorphisms, the co-unit $\vep : \cA\to \cD$ and the (oposit) co-multiplication $\Delta' : \cA\to \cA \widetilde{\otimes}\cA$ 
\bea
&&\vep(\hL^+_{i,j}(z))=\delta_{i,j}{T}_{Q_{\bep_i} }\quad (n\in \Z),
\qquad \vep(e^Q)=e^Q,\lb{counitUqp}\\
&&\vep(\mu_l({\hf}))= \vep(\mu_r(\hf))=\widehat{f}T_0, \lb{counitf}\\
&&\Delta'(\hL^+_{i,j}(z))=\sum_{k} 
\hL^+_{k,j}(z)\tot \hL^+_{i,k}(z),\lb{coproUqp}\\
&&\Delta'(e^{Q})=e^{Q}\tot e^{Q},\\
&&\Delta'(\mu_l(\hf))=1\widetilde{\otimes} \mu_l(\hf),\quad \Delta'(\mu_r(\hf))=\mu_r(\hf)\widetilde{\otimes} 1, \lb{coprof}
\ena
and  the algebra antihomomorphism $S:\cA\to \cA$ 
\be
&&S(\hL^+_{ij}(z))=(\hL^+(z)^{-1})_{ij},\\
&&S(e^{Q})=e^{-Q}, \quad S(\mu_r(\hat{f}))=\mu_l(\hat{f}),\quad S(\mu_l(\hat{f}))=\mu_r(\hat{f}).
\en

\section{Proof of Theorem \ref{actHC}}\lb{ProofacHC}
In this section we prove that the action of the half-currents \eqref{actK}-\eqref{actF} satisfies the defining relations of the elliptic algebra $E_{q,p}(\glnh)$ at $k=0$.  
From Lemma 6.10  in \cite{Konno16} it is enough to show the following relations, which are the dynamical counterpart  of those listed in 
Sec.C.1 in \cite{Konno16}  through Proposition \ref{dynamicalHC}.   
Or one can directly derive them from \eqref{DRLLpmpm} with corresponding Gauss decomposition of the dynamical $L$-operator.
 \bea
&&\cK_{j+1}^+(1/w_1)^{-1}E_{j+1,j}^+(1/w_2,P)\cK_{j+1}^+(1/w_1)\nn\\
&&=E_{j+1,j}^+(1/w_2,P+<\bep_{j+1},h>)\frac{1}{\bar{b}_{}^*(-v_{12})}
-E_{j+1,j}^+(1/w_1,P+<\bep_{j+1},h>)\frac{c_{}^*(-v_{12},P_{j,j+1})}{\bar{b}_{}^*(-v_{12})},\nn\\
&&\lb{dykek}\\
&&\cK_{j+1}^+(1/w_1)F_{j,j+1}^+(1/w_2,P+<\bep_{j+1},h>)\cK_{j+1}^+(1/w_1)^{-1}\nn\\
&&=\frac{1}{\bar{b}_{}(-v_{12})}F_{j,j+1}^+(1/w_2,P)-
\frac{\bar{c}_{}(-v_{12},(P+h)_{j,j+1})}{
\bar{b}_{}(-v_{12})}F_{j,j+1}^+(1/w_1,P), 
\lb{dykfk}\\
&&\frac{1}{\bar{b}^*(v_{12})}E^+_{j+1,j}(1/w_1,P+<\bep_{j+1},h>)E^+_{j+1,j}(1/w_2,P+<\bep_j,h>)\nn\\
&&\hspace{1cm}-E^+_{j+1,j}(1/w_2,P+<\bep_{j+1},h>)E^+_{j+1,j}(1/w_2,P+<\bep_j,h>)
\frac{{c}^*(v_{12},P_{j,j+1})}{\bar{b}^*(v_{12})}\nn\\
&&=\frac{1}{\bar{b}^*(-v_{12})}
E^+_{j+1,j}(1/w_2,P+<\bep_{j+1},h>)E^+_{j+1,j}(1/w_1,P+<\bep_j,h>)\nn\\
&&\hspace{1cm}
-E^+_{j+1,j}(1/w_1,P+<\bep_{j+1},h>)E^+_{j+1,j}(1/w_1,P+<\bep_j,h>)\frac{{c}^*(-v_{12},P_{j,j+1})}{\bar{b}^*(-v_{12})},
\lb{dyee}
\ena
\bea
&&\frac{1}{\bar{b}(-v_{12})}F^+_{j,j+1}(1/w_1,P)F^+_{j,j+1}(1/w_2,P)-F^+_{j,l}(1/w_1,P)^2
\frac{\bar{c}(-v_{12},(P+h)_{j,j+1}-2)}{\bar{b}(-v_{12})}\nn\\
&&=\frac{1}{\bar{b}(v_{12})}
F^+_{j,j+1}(1/w_2,P)F^+_{j,j+1}(1/w_1,P)-
F^+_{j,j+1}(1/w_2,P)^2\frac{\bar{c}(v_{12},(P+h)_{j,j+1}-2)}{\bar{b}(v_{12})},
\lb{dyff}\\
&&E^+_{j+1,j}(1/w_1,P)F^+_{j,j+1}(1/w_2,P+<\bep_j,h>)-F^+_{j,j+1}(1/w_2,P+<\bep_{j+1},h>)E^+_{j+1,j}(1/w_1,P)\nn\\
&&=\cK^+_{j}(1/w_2)\cK^+_{j+1}(1/w_2)^{-1}
\frac{\bar{c}^*(-v_{12},P_{j,j+1})}{\bar{b}^*(-v_{12})} \nn \\
&&\qquad\qquad\qquad\qquad\qquad-\cK^+_{j+1}(1/w_1)^{-1}\cK^+_{j}(1/w_1)
\frac{\bar{c}(-v_{12},(P+h)_{j,j+1})}{\bar{b}(-v_{12})},\lb{dyeffe}
 \ena
where $w_i=q^{2v_i}\ (i=1,2)$ and we set $v_{12}=v_1-v_2$. 
Note that at $k=0$, $p=p^*$, $r=r^*$, hence $b(u,P)=b^*(u,P), \bar{b}(u)=\bar{b}^*(u), c(u,P)=c^*(u,P), \bar{c}(u,P)=\bar{c}^*(u,P)$ etc.   

\noindent
{\it \eqref{dykek}}: Noting $<\bep_{j+1},h_{j,j+1}>=-1$, from \eqref{actK} and \eqref{actE} we have 
\be
&&{\rm LHS}=\sum_{k\in I_{j+1}}\frac{\bar{c}(u_k-v_2,P_{j,j+1})}{\bb(u_k-v_2)}\bb(u_k-v_1)^{-1}\prod_{l\in I_{j+1}\atop \not=k}\bb(u_{kl})^{-1}
\xi_{I^{k'}},
%\\
\en
\be
&&{\rm RHS}=\sum_{k\in I_{j+1}}\left(\frac{\bar{c}(u_k-v_2,P_{j,j+1}-1)}{\bb(u_k-v_2)}\frac{1}{\bb(-v_{12})}
-\frac{\bar{c}(u_k-v_1,P_{j,j+1}-1)}{\bb(u_k-v_1)}\frac{c(-v_{12},P_{j,j+1})}{\bb(-v_{12})}\right)\nn\\
&&\qquad\quad\qquad\times \prod_{l\in I_{j+1}\atop \not=k}\bb(u_{kl})^{-1}\xi_{I^{k'}}.
\en
Then the equality follows from the identity 
\bea
&&\frac{\bar{c}(-v_2,s)}{\bb(-v_2)}\bb(-v_1)^{-1}=
\frac{\bar{c}(-v_2,s-1)}{\bb(-v_2)}\frac{1}{\bb(-v_{12})}
-\frac{\bar{c}(-v_1,s-1)}{\bb(-v_1)}\frac{c(-v_{12},s)}{\bb(-v_{12})}.
\ena
\qed

One can prove \eqref{dykfk} similarly. 

\vspace{2mm}
\noindent
{\it \eqref{dyee}}: The action of the LHS on $\xi_I$ yields
\be
&&\sum_{a,b\in I_{j+1}\atop a\not=b} \left(
\frac{[v_{12}+1][P-1-u_b+v_1][P+1-u_a+v_2]}{[v_{12}][u_b-v_1][u_a-v_2]}\right.\nn\\
&&\left.-\frac{[P+v_{12}][1][P-1-u_b+v_2][P+1-u_a+v_2]}{[P][v_{12}][u_b-v_2][u_a-v_2]}
\right)\frac{[u_{ab}+1]}{u_{ab}}\prod_{k\in I_{j+1}\atop \not=a,b}\bb(u_{ak})^{-1}\bb(u_{bk})^{-1}\xi_{(I^{a'})^{b'}}\\
&&=-\frac{[P-1]}{[P]}\sum_{a,b\in I_{j+1}\atop a\not=b}
\frac{[P+1-u_a+v_2][v_1-u_b+P][v_2-u_b-1]}{[u_a-v_2][u_b-v_1][u_b-v_2]}\frac{[u_{ab}+1]}{[u_{ab}]}\nn\\
&&\qquad\qquad\qquad\times\prod_{k\in I_{j+1}\atop \not=a,b}\bb(u_{ak})^{-1}\bb(u_{bk})^{-1}\xi_{(I^{a'})^{b'}}
\en
where we set $P=P_{j,j+1}$ and $u_{ak}=u_a-u_k$ etc. The second equality follows from the identity
\be
&&[v_{12}+1][P-1-u_b+v_1][P][u_b-v_2]-[P+v_{12}][1][P-1-u_b+v_2][u_b-v_1]\nn\\
&&=-[P-1][v_{12}][v_1-u_b+P][v_2-u_b-1]. 
\en
Similarly  the RHS of \eqref{dyee} yields
\be
&&\sum_{a,b\in I_{j+1}\atop a\not=b} \left(
\frac{[v_{12}-1][P-1-u_b+v_2][P+1-u_a+v_1]}{[v_{12}][u_b-v_2][u_a-v_1]}\right.\nn\\
&&\left.-\frac{[P-v_{12}][1][P-1-u_b+v_1][P+1-u_a+v_1]}{[P][v_{12}][u_b-v_1][u_a-v_1]}
\right)\frac{[u_{ab}+1]}{u_{ab}}\prod_{k\in I_{j+1}\atop \not=a,b}\bb(u_{ak})^{-1}\bb(u_{bk})^{-1}\xi_{(I^{a'})^{b'}}\nn\\
&&=-\frac{[P-1]}{[P]}\sum_{a,b\in I_{j+1}\atop a\not=b}
\frac{[P+1-u_a+v_1][v_2-u_b+P][v_1-u_b-1]}{[u_a-v_1][u_b-v_1][u_b-v_2]}\frac{[u_{ab}+1]}{[u_{ab}]}\nn\\
&&\qquad\qquad\qquad\times\prod_{k\in I_{j+1}\atop \not=a,b}\bb(u_{ak})^{-1}\bb(u_{bk})^{-1}\xi_{(I^{a'})^{b'}}. 
\en
Taking the difference between the LHS and the RHS, we obtain
\bea
&&-\frac{[P-1]}{[P]}\sum_{a,b\in I_{j+1}\atop a\not=b}\left(
\frac{[P+1-u_a+v_2][v_1-u_b+P][v_2-u_b-1]}{[u_a-v_2][u_b-v_1][u_b-v_2]}\right.\nn\\
&&\left.\hspace{3cm}-\frac{[P+1-u_a+v_1][v_2-u_b+P][v_1-u_b-1]}{[u_a-v_1][u_b-v_1][u_b-v_2]}\right)\nn\\
&&\qquad\qquad\qquad\times\frac{[u_{ab}+1]}{[u_{ab}]}\prod_{k\in I_{j+1}\atop \not=a,b}\bb(u_{ak})^{-1}\bb(u_{bk})^{-1}\xi_{(I^{a'})^{b'}}\nn\\
&&=\frac{[P-1][P+1][v_{12}]}{[P]}\sum_{a,b\in I_{j+1}\atop a\not=b}\frac{[u_{ab}-1][u_{ab}+1]}{[u_{ab}]}f(u_a,u_b)\xi_{(I^{a'})^{b'}}\lb{LR}
\ena
where the equality follows from the identity
\be
&&\hspace{-1cm}-[P+1-u_a+v_2][v_1-u_b+P][v_2-u_b-1][u_a-v_1]+[P+1-u_a+v_1][v_2-u_b+P][v_1-u_b-1][u_a-v_2]\nn\\
&&=[P+1][v_{12}][u_{ab}-1][P+v_1+v_2-u_a-u_b]
\en
and we set
\be 
&&f(u_a,u_b)=\frac{[P+v_1+v_2-u_a-u_b]}{[u_a-v_1][u_a-v_2][u_b-v_1][u_b-v_2]}\prod_{k\in I_{j+1}\atop \not=a,b}\bb(u_{ak})^{-1}\bb(u_{bk})^{-1}.
\en
Since $f(u_a, u_b)=f(u_b,u_a)$ and $(I^{a'})^{b'}=(I^{b'})^{a'}\ (a,b\in I_{j+1}, a\not=b )$, 
the summation in \eqref{LR} vanishes. \qed

One can prove \eqref{dyff} similarly. 

\noindent
{\it \eqref{dyeffe}}: In the LHS, noting $<\bep_j,h_{j,j+1}>=1$ we obtain
\bea
&&E^+_{j+1,j}(1/w_1,P)F^+_{j,j+1}(1/w_2,P+<\bep_j,h>)\xi_I\nn\\
&&=\sum_{a\in I_j}\frac{c(u_a-v_2,P_{j,j+1}+\la_j-\la_{j+1})}{\bb(u_a-v_2)}\frac{\bar{c}(u_a-v_1,P_{j,j+1})}{\bb(u_a-v_1)}
\prod_{k\in I_j\atop \not=a}\bb(u_{ka})^{-1}\prod_{l\in I_{j+1}}\bb(u_{al})^{-1}\xi_{I}\nn\\
&&+\sum_{a\in I_j}\sum_{b\in I_{j+1}}\frac{c(u_a-v_2,P_{j,j+1}+\la_j-\la_{j+1})}{\bb(u_a-v_2)}\frac{\bar{c}(u_b-v_1,P_{j,j+1})}{\bb(u_b-v_1)}
\prod_{k\in I_j\atop \not=a}\bb(u_{ka})^{-1}\prod_{l\in I_{j+1}\atop \not=b}\bb(u_{bl})^{-1}\nn\\
&&\qquad\qquad\qquad \times \bb(u_{ba})^{-1}\xi_{(I^{'a})^{b'}},\lb{EFLHS}\\
&&F^+_{j,j+1}(1/w_2,P+<\bep_{j+1},h>)E^+_{j+1,j}(1/w_1,P)\xi_I\nn\\
&&=\sum_{b\in I_{j+1}}\frac{c(u_b-v_2,P_{j,j+1}+\la_j-\la_{j+1})}{\bb(u_b-v_2)}\frac{\bar{c}(u_b-v_1,P_{j,j+1})}{\bb(u_b-v_1)}
\prod_{l\in I_{j+1}\atop \not=b}\bb(u_{bl})^{-1}\prod_{k\in I_{j}}\bb(u_{kb})^{-1}\xi_{I}\nn\\
&&+\sum_{a\in I_j}\sum_{b\in I_{j+1}}\frac{c(u_a-v_2,P_{j,j+1}+\la_j-\la_{j+1})}{\bb(u_a-v_2)}\frac{\bar{c}(u_b-v_1,P_{j,j+1})}{\bb(u_b-v_1)}
\prod_{k\in I_j\atop \not=a}\bb(u_{ka})^{-1}\prod_{l\in I_{j+1}\atop \not=b}\bb(u_{bl})^{-1}\nn\\
&&\qquad\qquad\qquad\times \bb(u_{ba})^{-1}\xi_{(I^{b'})^{'a}},\lb{FELHS}
\ena
where we set $u_{ka}=u_k-u_a$ etc.. 
Since $(I^{'a})^{b'}=(I^{b'})^{'a}$ $\forall a\in I_j, \forall b\in I_{j+1}$, the second terms in \eqref{EFLHS} and \eqref{FELHS} coincide 
each other. Hence we have
\bea
&&\left(E^+_{j+1,j}(1/w_1,P)F^+_{j,j+1}(1/w_2,P+<\bep_j,h>)-
F^+_{j,j+1}(1/w_2,P+<\bep_{j+1},h>)E^+_{j+1,j}(1/w_1,P)\right)\xi_I\nn\\
&&=\left(\sum_{a\in I_j}\frac{c(u_a-v_2,P_{j,j+1}+\la_j-\la_{j+1})}{\bb(u_a-v_2)}\frac{\bar{c}(u_a-v_1,P_{j,j+1})}{\bb(u_a-v_1)}
\prod_{k\in I_j\atop \not=a}\bb(u_{ka})^{-1}\prod_{l\in I_{j+1}}\bb(u_{al})^{-1}\right.\nn\\
&&\left. -
\sum_{b\in I_{j+1}}\frac{c(u_b-v_2,P_{j,j+1}+\la_j-\la_{j+1})}{\bb(u_b-v_2)}\frac{\bar{c}(u_b-v_1,P_{j,j+1})}{\bb(u_b-v_1)}
\prod_{l\in I_{j+1}\atop \not=b}\bb(u_{bl})^{-1}\prod_{k\in I_{j}}\bb(u_{kb})^{-1}\right)\xi_{I}.\lb{EFFELHS}
\ena
In the RHS of \eqref{dyeffe}, noting ${\Delta'}^{(n-1)}(h_{j,j+1})\xi_I=(\la_j-\la_{j+1})\xi_{I}$  we obtain
\bea
&&\cK^+_j(1/w_2)\cK^+_{j+1}(1/w_2)^{-1}\frac{\bar{c}(-v_{12},P_{j,j+1})}{\bb(-v_{12})}\xi_{I}\nn\\
&&=\frac{\bar{c}(-v_{12},P_{j,j+1})}{\bb(-v_{12})}\prod_{a\in I_j}\bb(u_a-v_2)^{-1}\prod_{b\in I_{j+1}}\bb(v_2-u_b)^{-1}\xi_{I},
\lb{KjKjp1}\\
&&\cK^+_{j+1}(1/w_1)^{-1}\cK^+_{j}(1/w_1)\frac{\bar{c}(-v_{12},(P+{\Delta'}^{(n-1)}(h))_{j,j+1})}{\bb(-v_{12})}\xi_{I}\nn\\
&&= \frac{\bar{c}(-v_{12},P_{j,j+1}+\la_j-\la_{j+1})}{\bb(-v_{12})}\prod_{a\in I_j}\bb(u_a-v_1)^{-1}\prod_{b\in I_{j+1}}\bb(v_1-u_b)^{-1}\xi_{I}.
\lb{Kjp1Kj}
\ena
Hence 
\bea
{\rm LHS}-{\rm RHS}
&=&\left(\sum_{a\in I_j}\frac{c(u_{a,2},P_{j,j+1}+\la_j-\la_{j+1})}{\bb(u_{a,2})}\frac{\bar{c}(u_{a,1},P_{j,j+1})}{\bb(u_{a,1})}
\prod_{k\in I_j\atop \not=a}\bb(u_{ka})^{-1}\prod_{l\in I_{j+1}}\bb(u_{al})^{-1}\right.\nn\\
&&\qquad -
\sum_{b\in I_{j+1}}\frac{c(u_{b,2},P_{j,j+1}+\la_j-\la_{j+1}))}{\bb(u_{b,2})}\frac{\bar{c}(u_{b,1},P_{j,j+1})}{\bb(u_{b,1})}
\prod_{k\in I_{j}}\bb(u_{kb})^{-1}\prod_{l\in I_{j+1}\atop \not=b}\bb(u_{bl})^{-1}\nn\\
&&\qquad-\frac{\bar{c}(-v_{12},P_{j,j+1})}{\bb(-v_{12})}\prod_{a\in I_j}\bb(u_{a,2})^{-1}\prod_{b\in I_{j+1}}\bb(-u_{b,2})^{-1}\nn\\
&&\left.
\qquad+\frac{\bar{c}(-v_{12},P_{j,j+1}+\la_j-\la_{j+1}))}{\bb(-v_{12})}\prod_{a\in I_j}\bb(u_{a,1})^{-1}\prod_{b\in I_{j+1}}\bb(-u_{b,1})^{-1}
\right)\xi_{I},
\ena
where we set 
$u_{c,i}=u_c-v_i\ (c=a,b, i=1,2)$. 
For $a\in I_j$ let us regard $(\cdots)$ in the RHS as a function of $u_a$ and denote it by $F(u_a)$. 
It is not so hard to find 
\be
&&F(u_a+r)=F(u_a),\qquad F(u_a+r\tau)=e^{-\frac{2\pi i}{r}}F(u_a) 
\en
and all the residues at the poles $u_a=v_1, v_2, u_k, u_l \ (k\in I_j, k\not=a, l\in I_{j+1})$ vanish. 
Hence $F(u_a)$  should be identically zero, because  $F(u_a)[u_a]/[u_a+1]$ becomes a order 1 elliptic function 
unless $F(-1)=0$.  
\qed

\section{A Direct Check of  Corollary \ref{actUqp} for \eqref{EFFE} at $k=0$}\lb{proofEFFE}
In this section we give a direct check that the action of the elliptic currents in  Corollary \ref{actUqp} satisfies \eqref{EFFE}.  
We start from the following partial fraction expansion formula, which can be obtained from (4.2) in \cite{Rosengren} by
 changing the multiplicative notation to the additive one and setting $t=v$, $b_k=u_k-1\ (k=1,\cdots,m)$, $b_{l}=a_{l}+1\ (l=m+1,\cdots,n)$ and $b_{n+1}=v+2m-n$.  
\begin{lem}\lb{partialfrac}
\bea
&&\prod_{k=1}^m\frac{[v-u_k+1]}{[v-u_k]}\prod_{l=m+1}^n\frac{[v-u_l-1]}{[v-u_l]}\nn\\
&&=\sum_{a=1}^n\frac{[v-u_a+2m-n]}{[2m-n][v-u_a]}\frac{\prod_{k=1}^m[u_a-u_k+1]\prod_{l=m+1}^n[u_a-u_l-1]}
{\prod_{j=1\atop \not=a}^n [u_a-u_j] }. \lb{partialfrac}
\ena
\end{lem}

Let us denote the LHS of  \eqref{partialfrac} by $F(v)$. Then
\be
&&\mathop{\rm Res}_{v=u_a} F(v) d v=\frac{1}{[0]'}\frac{\prod_{k=1}^m[u_a-u_k+1]\prod_{l=m+1}^n[u_a-u_l-1]}{\prod_{j=1\atop \not=a}^n [u_a-u_j] }.
\en
Then from \eqref{pmdelta}, we obtain 
\begin{lem}\lb{Resformula}
\be
&&F(v)|_+-F(v)|_-=\sum_{a=1}^n\delta(z_a/w)\mathop{\rm Res}_{v=u_a} F(v)dv,
\en
where $z_a=q^{2u_a}, w=q^{2v}$.
\end{lem}

Now let us check the relation \eqref{EFFE}.   

\noindent
1) the cases $|i-j|>1$ and $i=j+1$: It is obvious that ${(I^{'a})^{b'}}={(I^{b'})^{'a}}$, $\forall a\in I_j$, $\forall b\in I_{i+1}$. Hence $[E_i(1/w_1),F_j(1/w_2)]\xi_I=0$. 

\noindent
2) the case $i=j-1$: It is easy to show ${(I^{'a})^{b'}}={(I^{b'})^{'a}}$, $\forall a, b\in I_j$, $(a\not=b)$. Then 
\be
&&E_{j-1}(q^{j-N}/w_1)F_j(q^{j-N+1}/w_2)\xi_I=C e^{Q_{\al_{j-1}}}\sum_{a,b\in I_j\atop a\not=b}\delta(z_a/w_2)\delta(z_b/w_1)
\prod_{k\in I_j\atop \not=a}\bb(u_{ka})^{-1}\prod_{l\in (I^{'a})_{j}\atop \not=b}\bb(u_{bl})^{-1}\xi_{(I^{'a})^{b'}},
\en
whereas
\be
&&F_j(q^{j-N+1}/w_2)E_{j-1}(q^{j-N}/w_1)\xi_I=C e^{Q_{\al_{j-1}}}\sum_{a,b\in I_j\atop a\not=b}\delta(z_a/w_2)\delta(z_b/w_1)
\prod_{k\in (I^{b'})_j\atop \not=a}\bb(u_{ka})^{-1}\prod_{l\in I_{j}\atop \not=b}\bb(u_{bl})^{-1}\xi_{(I^{b'})^{'a}}. 
\en
Here we set 
\bea
&&C=\mu\mu^*\left(\frac{[1]}{[0]'}\right)^2=-\frac{\varrho}{q-q^{-1}}\frac{[1]}{[0]'}.\lb{def:C}
\ena
The last equality follows from  \eqref{mumus}. Noting 
\be
&&\prod_{l\in I_{j}\atop \not=b}\bb(u_{bl})^{-1}=\bb(u_{ba})^{-1}\prod_{l\in  (I^{'a})_{j},
\atop \not=b}\bb(u_{bl})^{-1},
\en
we obtain $[E_{j-1}(q^{j-N}/w),F_j(q^{j-N+1}/w)]\xi_I=0$. 

\noindent
3) the case $i=j$:  We have
\bea
&&E_{j}(q^{j-N+1}/w_1)F_j(q^{j-N+1}/w_2)\xi_I\nn\\
&&=C e^{-Q_{\al_{j}}}\left(
\sum_{a\in I_j}\sum_{b\in I_{j+1}
}\delta(z_a/w_2)\delta(z_b/w_1)
\prod_{k\in I_j\atop \not=a}\bb(u_{ka})^{-1}\prod_{l\in I_{j+1}\atop \not=b}\bb(u_{bl})^{-1}\times \bb(u_{b,a})^{-1}\xi_{(I^{'a})^{b'}}\right.\nn\\
&&\left. 
\hspace{3cm}+\sum_{a\in I_j}\delta(z_a/w_2)\delta(z_a/w_1)
\prod_{k\in I_j\atop \not=a}\bb(u_{ka})^{-1}\prod_{l\in I_{j+1}}\bb(u_{al})^{-1}\xi_{I}\right),\lb{EFxi}
\ena
whereas
\bea
&&F_j(q^{j-N+1}/w_2)E_{j}(q^{j-N+1}/w_1)\xi_I\nn\\
&&=C e^{-Q_{\al_{j}}}\left(
\sum_{a\in I_j}\sum_{b\in I_{j+1}
}\delta(z_a/w_2)\delta(z_b/w_1)
\prod_{k\in I_j\atop \not=a}\bb(u_{ka})^{-1}\prod_{l\in I_{j+1}\atop \not=b}\bb(u_{bl})^{-1}\times \bb(u_{b,a})^{-1}\xi_{(I^{b'})^{'a}}\right.\nn\\
&&\left. 
\hspace{3cm}+\sum_{b\in I_{j+1}}\delta(z_b/w_2)\delta(z_b/w_1)
\prod_{k\in I_j}\bb(u_{kb})^{-1}\prod_{l\in I_{j+1} \atop \not=b}\bb(u_{al})^{-1}\xi_{I}\right).\lb{FExi}
\ena
Since $(I^{'a})^{b'}=(I^{b'})^{'a}$  $\forall a\in I_j, \forall b\in I_{j+1}$,  the first terms in \eqref{EFxi} and \eqref{FExi} coincide each other.
Hence we obtain
\bea
&&[E_{j}(q^{j-N+1}/w_1),F_j(q^{j-N+1}/w_2)]\xi_I\nn\\
&&=C e^{-Q_{\al_{j}}}\left(
\sum_{a\in I_j}\delta(z_a/w_2)\delta(z_a/w_1)
\prod_{k\in I_j\atop \not=a}\bb(u_{ka})^{-1}\prod_{l\in I_{j+1}}\bb(u_{al})^{-1}\right.\nn\\
&&\left.\hspace{3cm}
-\sum_{b\in I_{j+1}}\delta(z_b/w_2)\delta(z_b/w_1)
\prod_{k\in I_j}\bb(u_{kb})^{-1}\prod_{l\in I_{j+1} \atop \not=b}\bb(u_{al})^{-1}\right)
\xi_{I}.\lb{EFFExi}
\ena
Let us set the eigenvalue of $H^+_j(q^{j-N+1}/w)$ on $\xi_I$ by $\varrho h_j(v)$ with $\varrho$ in \eqref{def:kappa} i.e. 
\be
&& h_j(v)=\prod_{k\in I_j}\frac{[u_{k}-v+1]}{[u_k-v]}\prod_{l\in I_{j+1}}\frac{[u_{l}-v-1]}{[u_l-v]}.
\en
Then we have for $a\in I_j$ and $b\in I_{j+1}$, 
\be
&&\prod_{k\in I_j\atop \not=a}\bb(u_{ka})^{-1}\prod_{l\in I_{j+1}}\bb(u_{al})^{-1}=-\frac{[0]'}{[1]}\mathop{\rm Res}_{v=u_a}h_j(v)dv,\\
&&\prod_{k\in I_j}\bb(u_{kb})^{-1}\prod_{l\in I_{j+1} \atop \not=b}\bb(u_{al})^{-1}=\frac{[0]'}{[1]}\mathop{\rm Res}_{v=u_b}h_j(v)dv. 
\en
Therefore one can write \eqref{EFFExi} as
\bea
&&\hspace{-1cm}[E_{j}(q^{j}/w_1),F_j(q^j/w_2)]\xi_I
=-\frac{[0]'C}{[1]}e^{-Q_{\al_{j}}}\delta(w_1/w_2)\sum_{c\in I_j\cup I_{j+1}} \delta(z_c/w_1)
\mathop{\rm Res}_{v=u_c}h_j(v)dv\, \xi_I.\lb{EFFERes}
\ena
Then from Lemma \ref{Resformula} and \eqref{def:C},  we obtain the desired formula.
\be
[E_{j}(q^{j-N+1}/w_1),F_j(q^{j-N+1}/w_2)]\xi_I
&=&\frac{\varrho}{q-q^{-1}}\delta(w_1/w_2)\left(h_j(v)|_+-h_j(v)|_-\right)\xi_I\\
&=&\frac{1}{q-q^{-1}}\delta(w_1/w_2)\left(H^+_j(q^{j-N+1}/w_1)-H^-_j(q^{j-N+1}/w_1)\right)\xi_I.
\en
\qed

\end{appendix}

\renewcommand{\baselinestretch}{0.7}

\end{document}